\documentclass[final,hidelinks,onefignum,onetabnum]{siamart220329}

\usepackage{lipsum}
\usepackage{amsfonts}
\usepackage{graphicx}
\usepackage{epstopdf}
\usepackage{algorithmic}
\ifpdf
  \DeclareGraphicsExtensions{.eps,.pdf,.png,.jpg}
\else
  \DeclareGraphicsExtensions{.eps}
\fi

\newsiamremark{remark}{Remark}
\newsiamremark{hypothesis}{Hypothesis}
\crefname{hypothesis}{Hypothesis}{Hypotheses}
\newsiamthm{claim}{Claim}

\headers{A hybrid multiscale solver}{Michael Hintermüller and Denis Korolev}

\title{A hybrid physics-informed neural network based multiscale solver as a partial differential equation constrained optimization problem\thanks{Submitted to the editors DATE.
\funding{The authors acknowledge the support of the Leibniz Collaborative Excellence Cluster under project  ML4Sim (funding reference: K377/2021). MH also acknowledges support by the DFG ExC 2046 MATH+: Berlin Mathematics Research Center under project EF1-17.}}}

\author{Michael Hintermüller\footnotemark[3] \thanks{Institute for Mathematics, Humboldt-Universität zu Berlin, Unter den Linden 6, 10099 Berlin, Germany
  (\email{hint@mathematik.hu-berlin.de})} 
\and  Denis Korolev\thanks{Weierstrass Institute for Applied Analysis and Stochastics (WIAS), Mohrenstr. 39, 10117 Berlin, Germany
  (\email{korolev@wias-berlin.de, hintermueller@wias-berlin.de}).}
}

\usepackage{amsopn}

\ifpdf
\hypersetup{
  pdftitle={A hybrid physics-informed neural network based multiscale solver as a partial differential equation constrained optimization problem},
  pdfauthor={Michael Hintermüller, Denis Korolev}
}
\fi

\newsiamthm{Assumption}{Assumption} % <- Preamble
\newsiamthm{Example}{Example} % <- Preamble
\newsiamthm{Remark}{Remark} % <- Preamble
\newcommand\restr[2]{\ensuremath{\left.#1\right|_{#2}}}
\begin{document}

\maketitle

% REQUIRED
\begin{abstract}
In this work, we study physics-informed neural networks (PINNs)  constrained by partial differential equations (PDEs) and their application in approximating PDEs with two characteristic scales. From a continuous perspective, our formulation corresponds to a non-standard PDE-constrained optimization problem with a PINN-type objective. From a discrete standpoint, the formulation represents a hybrid numerical solver that utilizes both neural networks and finite elements. For the problem analysis, we introduce a proper function space, and we develop a numerical solution algorithm. The latter combines an adjoint-based technique for the efficient gradient computation with automatic differentiation. This new multiscale method is then applied exemplarily to a heat transfer problem with oscillating coefficients. In this context, the neural network approximates a fine-scale problem, and a coarse-scale problem constrains the associated learning process. We demonstrate that incorporating coarse-scale information into the neural network training process via a weak convergence-based regularization term is beneficial. Indeed, while preserving upscaling consistency, this term encourages non-trivial PINN solutions and also acts as a preconditioner for the low-frequency component of the fine-scale PDE, resulting in improved convergence properties of the PINN method.
\end{abstract}

% REQUIRED
\begin{keywords}
Learning-informed optimal control, PDE constrained optimization, physics-informed neural networks, quasi-minimization, weak convergence, multiscale modelling
\end{keywords}

% REQUIRED
\begin{MSCcodes}
35B27, 65K10, 65J15, 65N99, 68T20, 80M40 
\end{MSCcodes}

\section{Introduction}
Solving partial differential equations (PDEs) using physics-informed neural networks (PINNs) is currently an active area of research; see, e.g. \cite{cuomo2022scientific} for an overview and references therein. The main principle of physics-informed learning was pioneered by \cite{lagaris1998artificial} and later reincarnated in its modern computational interpretation by Raissi \textit{et al.} \cite{raissi2019physics}. It consists of integrating physical laws, typically in the form of the residuals of underlying PDEs, into a least-squares objective and finding the approximate solution to the corresponding residual minimization problem. The approximation ansatz $u_{\boldsymbol{\theta},n}$ for the PDE solution $u\in U$, with $U$ a suitable Banach space, is then sought in a neural network class $\mathfrak{N}_{\boldsymbol{\theta},n}$. The unknown network parameters $\boldsymbol{\theta} \in \mathbb{R}^{N_n}$, the so-called weights and biases, are determined by solving an associated nonlinear and non-convex optimization problem. PINN methods are typically meshless, rendering them potentially useful for numerically solving PDEs on complex domains or in high dimensions \cite{hu2023tackling, zang2020weak}. The associated framework is rather flexible and allows for easy data incorporation in various ways, making PINNs not only versatile, but also a  competitive approach for solving inverse problems \cite{chen2020physics, jagtap2022physics, mishra2022estimates}. In addition, the expressiveness of neural networks is supported by universal approximation theorems \cite{de2021approximation, guhring2020error, mhaskar1997neural, xie2011errors, yarotsky2017error} and transfer learning capabilities \cite{goswami2020transfer, xu2023transfer}. Moreover, (approximate) optimization can be performed rapidly on modern computers and compute clusters, thanks to the excellent parallelization capabilities of neural networks on GPUs, advances in automatic differentiation, as well as parallel and domain decomposition techniques \cite{jagtap2021extended,  moseley2023finite,  shukla2021parallel}.  On the other hand, the non-convex nature of the underlying optimization and complex nonlinear dynamics within the learning process often lead to difficulties and limitations rendering the analytical and numerical handling delicate \cite{wang2021understanding, wang2022and}. Furthermore, PINNs can be difficult to train for problems exhibiting high-frequency or multiscale behavior \cite{wang2021eigenvector}, particularly due to the so-called ``spectral bias'' of neural networks. The latter is related to the fact that learning prioritizes low-frequency modes and prevents networks from effectively learning high-frequency functions \cite{rahaman2019spectral}.  However, the combination of physics-informed neural networks with numerically robust and efficient solvers may yield a way to mitigate these challenges.

Motivated by multiscale systems, in this work we enhance neural network training by incorporating a learning-informed PDE as a constraint into the PINN optimization. For the ease of exposition, we consider here a two-scale setting only, which involves a fine-scale equation (of formidable computational complexity, perhaps beyond reach) and a coarse-scale equation (which we consider computationally tractable). The aforementioned computational burden may stem, e.g., from fine-scale properties of composite materials (foams, textiles, etc) resulting in highly oscillatory multiscale or high-contrast coefficients (heat conductivity, permeability, etc). A reliable associated simulation 
may then require a prohibitively fine numerical resolution. In order to remedy the enormous computational complexity, we use a neural network based approximation of the fine-scale problem which informs the coarse-scale problem through a homogenization procedure of choice. This gives rise to the following PDE-constrained optimization problem:
\begin{equation}\label{learning-informed problem}
\left\{\begin{split}
&{\inf} \ J(y, u_{\boldsymbol{\theta},n})\quad\text{over }(y , u_{\boldsymbol{\theta},n}), \\ 
&\text{subject to (s.t.) } \mathcal{L}\big[u_{\boldsymbol{\theta},n}\big] y= f, 
\end{split}\right.
\end{equation}
where $J$ stands for a (least-squares) loss functional penalizing the PDE residual of the fine-scale equation (possibly including boundary conditions). It is worth noting that in the standard PINN framework, $J$ depends solely on $u_{\boldsymbol{\theta},n}$. Here, however, we also introduce a coupling term to make the loss also dependent on $y$, thus stabilizing the training by the aforementioned coarse-scale enrichment. Conceptually, this additional term incorporates information on the weak convergence of the fine-scale solution to the coarse-scale one into the loss functional. By $\mathcal{L}[u_{\boldsymbol{\theta},n}]: Y \rightarrow Z$ we denote a coarse-scale differential operator between Banach spaces $Y$ and $Z$, which is informed by our neural network ansatz yielding $u_{\boldsymbol{\theta},n}$.  Together with some given data $f$, it defines an equality constraint in \eqref{learning-informed problem}. Structurally, $u_{\boldsymbol{\theta},n}$ may be considered a control variable, whereas the coarse-scale solution $y$ acts as a state.  The abstract framework is quite general, allowing for the use of various (multiscale) techniques to parameterize $\mathcal{L}\big[u_{\boldsymbol{\theta},n}\big]$ by $u_{\boldsymbol{\theta},n}$.

In the realm of learning-informed optimal control, several works, including \cite{dong2022optimization, sirignano2023pde}, have focused on approximating nonlinear constituents or source terms in the state equation using neural networks. PINNs have also been employed in solvers for underlying state and adjoint equations in various PDE-constrained optimization scenarios \cite{lu2021physics, mowlavi2023optimal}. In \cite{brevis2022neural}, the neural stabilization of non-stable discrete weak formulations is proposed, resulting in a non-standard PDE-constrained optimization with neural network controls.  Let us point out here that all the aforementioned techniques, while structurally perhaps similar, differ from this work. Indeed, in the usual approaches the objective is typically not related to a neural network learning problem or PDE residual minimization, but rather to minimizing specific (e.g., tracking-type) cost functionals \cite{casas1992optimal, hinze2008optimization}. This difference has significant implications in analysis and numerical implementation. Besides, to the best of the authors' knowledge, our work is the first to deal in detail with a PINN-based optimization problem constrained by a PDE. Here, problem \eqref{learning-informed problem} is formulated and analyzed in a function space setting, taking into account the regularity of the fine-scale and coarse-scale PDE solutions as well as the interdependent choices of activation functions and PINN losses. For our new problem \eqref{learning-informed problem}, the concept of quasi-minimization \cite{shin2023error} is crucial when studying (approximate) existence of solutions.

Discretizing the coarse-scale equation in \eqref{learning-informed problem}, e.g., via the finite element method, alongside a (meshless) PINN-based approach for the fine-scale problem, yields a {\it hybrid physics-informed multiscale} numerical solver. In this context, the meshless approach appears particularly useful for complex geometries. Note also that in the course of the optimization process for solving the hybrid finite dimensional approximate version of \eqref{learning-informed problem} possibly requires to frequently solve the discretized coarse-scale PDE. The hope now is that the coarse-scale equation can be solved numerically at a significantly lower cost (compute time) than computing the PINN solution, while still well informing low-frequency components of the fine-scale solution. Our numerical experiments provide evidence that incorporating the coarse-scale solution into the learning process through a coupling term in the objective of \eqref{learning-informed problem} acts as a preconditioner for the low-frequency component of the fine-scale PDE, thereby accelerating convergence to such a component and improving the overall performance. The proposed methodology can be used with most PINN architectures (standard PINNs \cite{raissi2019physics}, Fourier features networks \cite{tancik2020fourier, wang2021eigenvector}, FBPINNs \cite{moseley2023finite}, etc) and has the potential to improve existing benchmarks. 

PINNs find numerous applications in multiscale systems and material design (see, e.g., \cite{chen2020physics, diao2023solving, lu2021physics, zhang2022analyses}). Some NN-based homogenization techniques \cite{arbabi2020linking, han2023neural, kropfl2023neural}, including PINN-based approaches \cite{leung2022nh, park2022physics}, have also been proposed. Besides, it is worth noting that decades of research in multiscale modelling have yielded powerful methods for obtaining surrogates at the coarse scale (see, e.g., \cite{abdulle2012heterogeneous, efendiev2013generalized, iliev2011variational, maalqvist2014localization, weinan2003heterogeneous}). In our application section, we utilize the upscaling method \cite{farmer2002upscaling, wu2002analysis} for numerical homogenization, facilitating the transition from fine-scale to coarse-scale domains via our parametrization $\mathcal{L}\big[u_{\boldsymbol{\theta},n}\big]$. However, despite sharing a similar application context with the aforementioned works, our focus here lies in enhancing PINN approximation for fine-scale problems by integrating it with numerical coarse-scale solvers using multiscale modelling techniques. The resulting hybrid coarse-coarse scale approximation is a byproduct of the neural network approximation on a fine scale and the chosen multiscale coupling, leading to the related \textit{upscaling consistency} concept.

The paper is organized as follows. In Section~\ref{Section 2}, we introduce an abstract framework for our learning-informed optimization problem. We recall specific details on physics-informed neural networks and introduce a compression operator, which embeds information on weak convergence into our loss function. Then, we derive the related upscaling consistency result. In Section~\ref{Section 3}, we fit \eqref{learning-informed problem} into a function space framework and propose a numerical algorithm for its solution. Section~\ref{Section 4} presents an upscaling technique for numerical homogenization. We then integrate the upscaling lift into the learning-informed PDE-constrained optimization setting and apply it to the heat transfer problem with oscillating coefficients.

% The outline is not required, but we show an example here.

\section{A hybrid multiscale approach}
\label{Section 2}
In this section, we define a function space and a coupling framework for two PDEs, the fine-scale and the coarse-scale problem, respectively. For the treatment of the associated PINNs, we closely follow \cite{shin2023error}. For the bounded domain $\Omega \subset \mathbb{R}^{d}$ with Lipschitz boundary $\partial \Omega$, let $L^{p}(\Omega)$, $H^{1}(\Omega)$, $H^{1}_{0}(\Omega)$, $H^{k}(\Omega)$, $W^{k,p}(\Omega)$, etc. denote the standard Lebesgue and Sobolev spaces; see, e.g., \cite{adams2003sobolev}. We also set $\mathbb{R}_{+}:= \{x \in \mathbb{R}: \ x > 0\}$ and $\mathbb{R}_{\geq 0}:= \{x \in \mathbb{R}: \ x \geq 0\}$.

\subsection{Function spaces and PDEs} Let $(U, \lVert \cdot \lVert_{U})$, $(H, \lVert \cdot \rVert_{H})$ be Hilbert spaces, $(X, \lVert \cdot \rVert_{X})$ be a normed vector space, $(V, \lVert \cdot \rVert_{V})$, $(Z, \lVert \cdot \rVert_{Z})$ be Banach spaces, and $X$ a dense subspace of $U$. Suppose also that $U$ is continuously embedded in $V$ and $V$ is continuously embedded in $H$, i.e., $U \hookrightarrow V \hookrightarrow H$. Moreover, $X \hookrightarrow U$. 

For given $f^{\varepsilon} \in H$, we consider the following partial differential equation
\begin{align}\label{MP: eq 1}
\mathcal{A}^{\varepsilon}u= f^{\varepsilon}, \quad \text{in} \   H, \quad
\mathcal{B}u = 0, \quad \text{in} \  Z, 
\end{align}
with $\mathcal{A}^{\varepsilon}: U \rightarrow H$ a bounded linear partial differential operator, i.e., $\mathcal{A}^{\varepsilon}\in L(U,H)$, depending on the fine-scale length $\varepsilon>0$, and $\mathcal{B}\in L(U,Z)$ defining boundary conditions. For \eqref{MP: eq 1} we invoke the following solution concept.
\begin{Assumption} 
\label{Ass:existence}
\textit{For every $\varepsilon>0$, \eqref{MP: eq 1} admits a unique solution $u^{\varepsilon} \in U$ for which there exists an approximating sequence $\{u_{k}^{\varepsilon}\} \subset X$} such that 
\begin{align}\label{MP: eq 3}
    \underset{k \rightarrow \infty}{\lim} \ \lVert u_{k}^{\varepsilon} - u^{\varepsilon}  \rVert_{U} = 0, \ \ \underset{k \rightarrow \infty}{\lim} \ \lVert \mathcal{A}^{\varepsilon} u_{k}^{\varepsilon} - f^{\varepsilon}  \rVert_{H} + \lVert \mathcal{B} u_{k}^{\varepsilon}  \rVert_{Z}  = 0.
\end{align}
%where $\{u_{k}^{\varepsilon}\} \subset X$ is an approximating sequence of $u^{\varepsilon}$.
\end{Assumption}
\begin{Remark} \label{Non-homogeneous bc remark} Let $g\in Z$ and assume that there exists $u^g\in U$ such that $\mathcal{B}u^g=g$ in $Z$. If $u^\varepsilon\in U$ satisfies \eqref{MP: eq 3} with $f^{\varepsilon}=q-\mathcal{A}^{\varepsilon}u^g$ for some given $q\in H$, then $w^{\varepsilon}:=u^\varepsilon + u^g\in U$ corresponds to \eqref{MP: eq 1} with $\mathcal{A}^{\varepsilon}u= q$ in $H$ and $\mathcal B u=g$ in $Z$.
\end{Remark}
\noindent
From now on we assume that Assumption~\ref{Ass:existence} is satisfied.
%, whenever this is needed. 
Note that if $u^{\varepsilon} \in U$ is the solution to \eqref{MP: eq 1}, then $u^{\varepsilon} \in U_0: = \{ u \in U: \ \mathcal{B}u = 0 \} \subseteq U$, and we sometimes use $\lVert \cdot \rVert_{U_0}$ instead of $\lVert \cdot \rVert_{U}$. 
%The norm $\lVert \cdot \rVert_{\widetilde{U}}$ on $\widetilde{U}$ coincides with the norm on $U$, and hence they are used interchangeably. 
The following stability estimate is standard for least-squares residual minimization, including the least-squares finite element method \cite{bochev2016least, bochev2009least} and physics-informed neural networks \cite{mishra2022estimatesforward, shin2023error, zeinhofer2023unified}.    
\begin{Assumption}[\textbf{Fine-scale stability}] 
\label{Ass:Stability}
\textit{There exist a stability bound $C_{s}^\varepsilon \in \mathbb{R}_{+}$ and an upper bound $C_{b}^\varepsilon \in \mathbb{R}_{+}$, both possibly dependent on $\varepsilon>0$, such that
%$\mathcal{A}^{\varepsilon}$ and $\mathcal{B}$:
\begin{align}\label{MP: subeq 2(a)}
C_{s}^\varepsilon \lVert u \rVert_{V}^{2} \leq \lVert \mathcal{A}^{\varepsilon}u \rVert_{H}^{2} + \lVert \mathcal{B} u \rVert_{Z}^{2} \leq C_{b}^\varepsilon \lVert u \rVert_{U}^{2}, \quad \forall{u} \in U. 
\end{align}
%where $C_s^\varepsilon$ and $C_b^\varepsilon$ may depend on $\varepsilon>0$.
}
\end{Assumption}
\noindent
The stability bound \cref{MP: subeq 2(a)} is well-suited for PINN problems which include boundary conditions via penalization. In this case a term penalizing violations of the boundary conditions is added to the PINN objective. %boundary conditions, e.g., by including boundary residual penalty terms in the least-squares loss. 
Note, however, that boundary conditions may also be imposed exactly, i.e., as (hard) constraints, for PINNS \cite{lagari2020systematic, lu2021physics, sukumar2022exact}. Then, one may ask \eqref{MP: subeq 2(a)} to hold for all $u \in U_{0}$, while dropping the term $\lVert \mathcal{B} u \rVert_{Z}^{2}$ and using a stronger norm $\lVert \cdot \rVert_{U}$ in the lower bound. We demonstrate both approaches in our example section and use \eqref{MP: subeq 2(a)} for our abstract formulation. 

Let $(Y, \lVert \cdot \rVert_{Y})$ be a Hilbert space with $U_{0} \subseteq Y$, $Y^{*}$ the topological dual space of $Y$ and $Y \hookrightarrow H\cong H^* \hookrightarrow Y^{*}$ a Gelfand triple with the compact embedding $Y \hookrightarrow H$. Let $\mathcal{L}[u]\in L(Y,Y^{*})$ be parameterized by $u \in U$. The corresponding bilinear form $b_{\mathcal{L}}[u] : Y \times Y \rightarrow \mathbb{R}$ is defined by
\begin{align}\label{MP: eq 3 forms}
b_{\mathcal{L}}[u](v,w):= \langle \mathcal{L}[u]v, w \rangle_{Y^{*}, Y}
\end{align}
for $v,w\in Y$. The following assumption is needed to secure (local) solvability of the coarse-scale equation.
\begin{Assumption}[\textbf{Coarse-scale stability}]
\label{Ass:Uniformity}
Let $u^{\varepsilon} \in U$ be the solution of \eqref{MP: eq 1}, $\bar{r}^\varepsilon \in \mathbb{R}_{+}$ and $B_{\bar{r}^\varepsilon}(u^{\varepsilon}):=\{v \in U: \ \lVert u^{\varepsilon} - v \rVert_{U} \leq \bar{r}^\varepsilon\}$. We assume that the form \eqref{MP: eq 3 forms} is uniformly bounded and uniformly coercive for all $u \in B_{\bar{r}^{\varepsilon}}(u^{\varepsilon})$, i.e., there exist $C_{b}^{\varepsilon}, C_{c}^{\varepsilon}\in\mathbb{R}_+$, independent of $u$, but possibly dependent on $\varepsilon$, such that  
\begin{align}\label{MP: eq 3.1 forms}
b_{\mathcal{L}}[u](v,w) \leq C_{b}^{\varepsilon} \lVert v \rVert_{Y} \lVert w \rVert_{Y}, \ \text{and} \ \  b_{\mathcal{L}}[u](v,v) \geq C_{c}^{\varepsilon} \lVert v \rVert^{2}_{Y}, \quad \forall{v,w}\in Y.
\end{align}
\end{Assumption}
\noindent
For given $(u^{\varepsilon}, \tilde{f}) \in U \times H$, we consider the following partial differential equation, which we refer to as the coarse-scale problem: Find $y(u^{\varepsilon}) \in Y$ such that
\begin{align}\label{MP: eq 4}
b_{\mathcal{L}}[u^{\varepsilon}](y(u^{\varepsilon}),v) = \langle \tilde{f}, v \rangle_{Y^{*}, Y} \quad \forall{v} \in Y. 
\end{align}
The coarse-scale problem \eqref{MP: eq 4} is well-posed by the Lax--Milgram lemma. Indeed, there exist a unique solution $y(u^{\varepsilon})\in Y$ and a bound $C^{\varepsilon} \in \mathbb{R}_{+}$ such that $\lVert y(u^{\varepsilon}) \rVert_{Y} \leq C^{\varepsilon} \lVert \tilde{f} \Vert_{Y^{*}}$, with $C^{\varepsilon}$ possibly dependent on $\varepsilon$ due to \eqref{MP: eq 3.1 forms}. 

\subsection{Physics-informed neural networks} For $L\in\mathbb{N}$, an $L$-layer feed-forward neural network (NN) is a recursively defined function $R_{\boldsymbol{\theta}}(x): \mathbb{R}^{n_{0}} \rightarrow  \mathbb{R}^{n_{L}} $ with
\begin{align*}
R_{\boldsymbol{\theta}}(x) = z^{L}(x), \ z^{l}(x)=W^{l}\sigma (z^{l-1}(x)) + b^{l},  \ 2 \leq l \leq L, \ z^{1}(x) = W^{1}x + b^{1},
\end{align*}
where $W^{l} \in \mathbb{R}^{n_{l} \times n_{l-1}}$ is the $l$-th layer weight matrix, $b_{l} \in \mathbb{R}^{n_{l}}$ is the $l$-th layer bias vector and $\sigma:\mathbb{R}\to\mathbb{R}$ is the activation function, which is applied component-wise in case of input arguments in $\mathbb{R}^{n_l}$. The network architecture is represented by the vector $\vec{\boldsymbol{n}}= (n_{0},...,n_{L})$, and the set of all possible network parameters is defined by
\begin{align*}
    \Theta(\vec{\boldsymbol{n}}) = \big \{ \{(W_{j}, b_{j})\}_{j=1}^{L} : W_{j} \in \mathbb{R}^{n_{j} \times n_{j-1}}, \  b_{j}  \in \mathbb{R}^{n_{j}} \big \}.
\end{align*}
For a sequence $\{\vec{\boldsymbol{n}}_{n} \}_{n\geq 1}$ of network architectures such that $\vec{\boldsymbol{n}}_{n} \leq \vec{\boldsymbol{n}}_{n+1}$ for all $n$, with the vector inequality understood component-wise, we define its corresponding sequence of neural network classes \cite{shin2023error}, and the related maximal realization map
\begin{align}\label{parameter iso}
\mathfrak{N}_{\theta, n}: = \{ R_{\boldsymbol{\theta}}: \boldsymbol{\theta} \in \cup_{\vec{\boldsymbol{v}} \leq \vec{\boldsymbol{n}}_{n}} \Theta(\vec{\boldsymbol{v}})\},  \ \text{and} \ \ \mathfrak{F}_{n} : \mathbb{R}^{N_{n}} \rightarrow  \mathfrak{N}_{\theta, n}, \quad \boldsymbol{\theta} \mapsto \mathfrak{F}_{n}(\boldsymbol{\theta}) =:v_{\boldsymbol{\theta},n}, 
\end{align} 
where $N_{n}$ denotes the total number of parameters in $\Theta(\vec{\boldsymbol{n}}_{n})$. Note that $\mathfrak{N}_{\theta, n} \subset \mathfrak{N}_{\theta, n+1}$ by construction, and that the regularity of $\mathfrak{F}_{n}$ depends on the one of the underlying activation functions in $\mathfrak{N}_{\theta, n}$. In what follows, $v_{\boldsymbol{\theta},n} \in \mathfrak{N}_{\theta, n}$ refers to a network function, generated by its maximal admissible number of parameters $N_{n}$.

The residual minimization of the fine-scale PDE over the class of neural networks leads to the following PINN optimization problem: 
\begin{align} \label{Pinn_opt}
%\underset{v_{\boldsymbol{\theta},n} \in \mathfrak{N}_{\theta, n} \cap U}{
\inf \mathcal{J}(v_{\boldsymbol{\theta},n}):= \lVert \mathcal{A}^{\varepsilon} v_{\boldsymbol{\theta},n} - f^{\varepsilon} \rVert_{H}^{2} +  \tau_{1}\lVert \mathcal{B}v_{\boldsymbol{\theta},n}\rVert_{Z}^{2}\quad\text{over }v_{\boldsymbol{\theta},n} \in \mathfrak{N}_{\theta, n},
\end{align}
where the penalty parameter $\tau_1 >0$ is fixed. We observe that $\mathcal{J}(u) \geq 0$ for all $u \in U$, and  $\mathcal{J}(u^{\varepsilon})=0$ for the solution $u^{\varepsilon}$ to \eqref{MP: eq 1}. Assuming momentarily that $U=H^{2}(\Omega)$ and $\mathcal{A}^{\varepsilon}$ is of second order, we need to ensure that $u^{\varepsilon}$ can be approximated by a neural network. It is well-known that Sobolev functions can be well approximated via popular deep ReLU neural networks; see, e.g., \cite{guhring2020error}.  Note, however, that the ReLU activation function $\sigma(x) = \max (0,x)$ admits only one weak derivative, which renders it infeasible when utilizing the standard least-squares loss with $H=L^{2}(\Omega)$ and $\mathcal{A}^{\varepsilon}$ involving higher-order (weak) derivatives. The hyperbolic tangent function $\tanh(x)$ can then be used to approximate functions from Sobolev spaces $H^{k}(\Omega)$, $k \geq 3$, in the norm of $U$ (see \cite[Theorem B.7]{de2022error}, but also \cite{de2021approximation}). We mention that the regularity requirements on activation functions and $u^{\varepsilon}$ can be relaxed by adopting a variational loss function \cite{kharazmi2019variational, kharazmi2021hp}. This involves multiplying the residuals of \eqref{MP: eq 1} by suitably smooth test functions and integrating by parts. Then, the respective weak residual is minimized over a (discrete) trial class of neural networks. In such a setting, one is confronted with the additional burden of discretizing the space of test functions. Typically, piecewise polynomials of low order, leading to a Petrov-Galerkin discretization, are preferred \cite{berrone2022variational}. Moreover, the Fourier transform can be employed to efficiently evaluate the $H^{-1}$ norm of the PDE residual for simple geometries \cite{taylor2023deep}.

In view of the above, in our work we focus on a setting with smooth activation functions. Instead of imposing assumptions on the sufficiently high number of admissible weak derivatives of $u^{\varepsilon}$, we rather work in a dense subspace $X$ of $U$, where the approximation with neural networks is less restrictive, such as $X=C^{k}(\bar{\Omega})$, which motivates the triplet of spaces $\mathfrak{N}_{\theta, n} \subset X \subset U$ and our definition of solution \eqref{MP: eq 3}. Then, we are interested in neural networks that can well approximate elements in $X$; cf. \cite{shin2023error} and some universal approximation results \cite{de2021approximation, guhring2020error, mhaskar1997neural, xie2011errors, yarotsky2017error}.
\begin{Assumption}[\textbf{Uniform NN approximation of elements in $X$}]
\label{Ass:Uniform NN Approximation of elements}
\textit{There exists a sequence of neural network classes $\{\mathfrak{N}_{\theta, n}\}$ with $\mathfrak{N}_{\theta, n}\subset X$ and $\mathfrak{N}_{\theta, n} \subset \mathfrak{N}_{\theta, n+1}$ for all $n\in\mathbb{N}$, and $X\subset\overline{\cup_{n}\mathfrak{N}_{\theta, n}}$ in the topology of $(X, \lVert \ \cdot \  \rVert_{X})$. In addition, it holds that $\mathcal{A}^{\varepsilon}v_{\boldsymbol{\theta}, n} \in L^{2}(\Omega)$ and $\mathcal{B}v_{\boldsymbol{\theta}, n} \in L^{2}(\partial \Omega)$ for all $v_{\boldsymbol{\theta}, n} \in \mathfrak{N}_{\theta, n}$.  %\textcolor{blue}{needs discussion}
}
\end{Assumption}

\subsection{Weak convergence based regularization} Let $H=L^{2}(\Omega)$ and $H^{\delta}:= L^{2}(\Omega_{\delta})$, where $\Omega_\delta:= \{z \in \Omega : \textnormal{dist}(z, \partial \Omega) > \frac{\delta}{2} \} \subset \Omega$  for $\delta > 0$.  The weak convergence assumption is common in the context of homogenization, and we make use of the following.  
\begin{Assumption}[\textbf{Weak convergence}] 
\label{Ass: Weak convergence}
\textit{Assume that for $\gamma > 0$, $u^{\varepsilon} \rightharpoonup y^{0}$ in $L^{2+\gamma}(\Omega)$ as $\varepsilon \rightarrow 0$, where  $u^{\varepsilon}$ is the solution to \eqref{MP: eq 1} and $y^{0} \in Y$ is the solution to the homogenized equation $\mathcal{L}^{0}y^{0}= f^{0}$, with $\mathcal{L}^{0} \in L(Y,Y^{\ast})$ and $f^{0} \in H$ $\varepsilon$-independent, respectively.}
\end{Assumption}
\begin{Remark} In our setting, $f^{\varepsilon}$ in \eqref{MP: eq 1}, $\tilde{f}$ in \eqref{MP: eq 4}, and $f^{0}$ in Assumption~\ref{Ass: Weak convergence} are different in general. This is due to the lifting of Dirichlet boundary conditions in all the PDEs, cf. Remark~\ref{Non-homogeneous bc remark}. Besides, we assume that $y^{0}$ has at least the same regularity as $y(u^{\varepsilon})$, as it follows from Assumption \ref{Ass: Weak convergence} and our choice of spaces. 
\end{Remark}
\noindent
In addition, we need the following; cf, e.g. \cite[Theorem 4.5]{abdulle2012heterogeneous} or \cite[Section 3.2]{wu2002analysis}.
\begin{Assumption}[\textbf{Consistent parametrization}] 
\label{Ass: Consistent parametrization}
\textit{There exist $C_{\mathcal{L}[u]}, C_\varepsilon \in \mathbb{R}_{+}$ such that $\lVert y^{0}  -  y(u^{\varepsilon})\rVert_{H} \leq C_{\mathcal{L}[u]} + C_\varepsilon$, where $C_{\mathcal{L}[u]}$ does not depend on $\varepsilon$ and can be made small by an appropriate parametrization of $y(u^{\varepsilon})$ by $u^{\varepsilon}$, and $C_\varepsilon \rightarrow 0$ as $\varepsilon \rightarrow 0$.}
\end{Assumption}
%Set $Q_{\delta} : H \rightarrow H^{\delta}$ as follows: 
For $v \in H$, $x \in \Omega_{\delta}$ and $\mathcal{V}_{\delta}(x) = \{z : \lVert z - x \rVert_{\mathbb{R}^{d}}  \leq \frac{\delta}{2} \} \subset \Omega$, define $Q_{\delta} : H \rightarrow H^{\delta}$ by
\begin{align*}
(Q_{\delta}v)(x) :=  \frac{1}{|\mathcal{V}_{\delta}(x)|} \int_{\mathcal{V}_{\delta}(x)} v(z) \  dz. 
\end{align*}
Following heterogeneous multiscale methods \cite{abdulle2012heterogeneous}, we fix the compression operator
\begin{align}\label{Compression operator}
(\bar{Q}_{\delta}v)(x) := 
\begin{cases}
                                   (Q_{\delta}v)(x), &\   x\in \Omega_{\delta}, \\
                                   v(x), &\  x \in \Omega \setminus  \Omega_{\delta}.
\end{cases}
\end{align}
Let us next introduce the coupling term $\mathcal{R}_{\delta}: Y \times U \rightarrow \mathbb{R}_{\geq 0}$ as follows:  
\begin{align}\label{MP: eq 8}
\mathcal{R}_{\delta}(y(u^{\varepsilon}), u^{\varepsilon}): = \lVert \bar{Q}_{\delta}u^{\varepsilon} - y(u^{\varepsilon}) \rVert_{H}^{2}.
\end{align}
The purpose of \eqref{MP: eq 8} is to equip the optimization problem \eqref{Pinn_opt} with information on weak convergence of the fine-scale solution to the coarse-scale one. 

Next we study \eqref{MP: eq 8}. We start with two preparatory result.
\begin{lemma}
\label{L2 convergence lemma}
The operator $Q_{\delta}:H \rightarrow H^{\delta}$ has the following properties: 
\begin{enumerate}
    \item $Q_{\delta} \in L(H, H^{\delta})$ and $\bar{Q}_{\delta} \in L(H):=L(H,H)$.
    \item Suppose that $u^{\varepsilon} \rightharpoonup y^{0}$ in $H$ as $\varepsilon \rightarrow 0$. Then, $\underset{\varepsilon \rightarrow 0}{\lim }\ \lVert Q_{\delta}u^{\varepsilon} - Q_{\delta}y^{0} \rVert_{H^{\delta}} = 0$.  
\end{enumerate}
\end{lemma}
\begin{proof}
The linearity of $Q_{\delta}$ is obvious. For $x \in \Omega_{\delta}$, using the Cauchy--Schwarz inequality, we obtain the estimate $|(Q_{\delta}u^{\varepsilon})(x)| \leq |\mathcal{V}_{\delta}(x)|^{-1/2} \lVert u^{\varepsilon}\rVert_{H}$. Since $|\mathcal{V}_{\delta}(x)|$ is constant for all $x \in \Omega_{\delta}$, the estimate is uniform.  By integrating the estimate over $\Omega_{\delta}$, we get $\lVert Q_{\delta}u^{\varepsilon} \rVert_{H^{\delta}}  \leq C \lVert u^{\varepsilon} \rVert_{H}$, where $C \in \mathbb{R}_{+}$. Therefore, $Q_{\delta} \in L(H, H^{\delta})$ and we readily establish that $\bar{Q}_{\delta} \in L(H)$.

Let $u^{\varepsilon} \rightharpoonup y^{0}$ in $H$ as $\varepsilon \rightarrow 0$. Then $\lVert u^{\varepsilon}\rVert_{H} \leq C$, where $C\in \mathbb{R}_{+}$ is $\varepsilon$-independent, and  for all test functions $v \in H$ it holds that
\begin{align}\label{MP: eq 5}
\int_{\Omega} u^{\varepsilon}(x) v(x) \  dx  \rightarrow \int_{\Omega} y^{0}(x) v(x)  \  dx  \quad \text{as} \quad  \varepsilon \rightarrow 0.
\end{align}
Consider the normalized characteristic function 
\begin{align*}
\bar{\chi}_{\mathcal{V}_{\delta}(x)}(z) := 
\begin{cases}
                                   \frac{1}{|\mathcal{V}_{\delta}(x)|}, &\   z \in \mathcal{V}_{\delta}(x), \\
                                   0,  &\  z \notin \mathcal{V}_{\delta}(x),
\end{cases}
\end{align*}
as the test function in \eqref{MP: eq 5}. Then one obtains the pointwise convergence of averages $(Q_{\delta}u^{\varepsilon})(x) \rightarrow (Q_{\delta}y^{0})(x)$ as $\varepsilon \rightarrow 0$ for $x\in \Omega_{\delta}$. The $L^{2}$ - convergence follows from the uniform estimate of $|(Q_{\delta}u^{\varepsilon})(x)|$ and Lebesgue's Dominated Convergence Theorem.
\end{proof}
\noindent
Since $Y \hookrightarrow H$, for a given $y \in Y$, almost every $x \in \Omega_{\delta}$ is a Lebesgue point, i.e.,
\begin{align}\label{MP: eq 7.1}
\underset{\delta \rightarrow 0}{\lim} \  \frac{1}{|\mathcal{V}_{\delta}(x)|}\int_{\mathcal{V}_{\delta}(x)}y(z) \ dz = y(x).
\end{align}
For small $\delta>0$, we consider the following approximation of the above limit: 
\begin{align*}
\frac{1}{|\mathcal{V}_{\delta}(x)|}\int_{\mathcal{V}_{\delta}(x)}y(z) \ dz \approx \underset{\delta \rightarrow 0}{\lim} \  \frac{1}{|\mathcal{V}_{\delta}(x)|}\int_{\mathcal{V}_\delta(x)}y(z) \ dz  = y(x).
\end{align*}
Additional integrability of $\nabla y$ provides a rate for such an approximation.
\begin{lemma}
\label{Lebesque approximation decay lemma}
Suppose that $y \in W^{1,p}(\Omega)$ with $\Omega \subset \mathbb{R}^{d}$ and $d < p < \infty$. Then
\begin{align*}
\lVert  y - Q_{\delta} y \rVert_{H^{\delta}}  \leq  C_{p}(y) \delta^{\frac{d(p-2)+1}{p}},
\end{align*}
where the constant $C_{p}(y) < \infty$ is independent of $\delta$, but it depends on $\lVert \nabla y \rVert_{L^{p}(\Omega)}$. 
\end{lemma}
\begin{proof}
First, we use a well-known trick that controls the deviation of a function from its average on convex sets (see e.g. \cite[Lemma 4.28]{adams2003sobolev}). Let $x \in \Omega_{\delta}$ and define the convex ball $\mathcal{V}_{\delta}(x) \subset \Omega$. For $y \in C^{1}(\bar{\Omega})$, $z \in \mathcal{V}_{\delta}(x)$ and for all $t\in [0,1]$, we get
\begin{align*}
y(z) - y(x) = \int_{0}^{1} \frac{d}{dt} y (x + t (z-x))  dt = \int_{0}^{1} \nabla y (x + t(z-x))\cdot (z-x)  dt.
\end{align*}
Integrating the above equality over $\mathcal{V}_{\delta}(x)$ and performing the change of variables $\xi=x + t(z-x)$, we obtain
\begin{align}\label{MP: eq 7.2}
\Big|\int_{\mathcal{V}_{\delta}(x)} y(z) dz - |\mathcal{V}_{\delta}(x)|\  y(x) \Big| \leq \int_{0}^{1} \frac{1}{t^{d}} \int_{\mathcal{V}_{t\delta}(x)} |\nabla y (\xi)| |\frac{\xi - x}{t}| \ d\xi\, dt. 
\end{align}
Hölder's inequality with $\frac{1}{p}+\frac{1}{q}=1$, in conjunction with the fact that $\mathcal{V}_{t \delta}(x) \subseteq \mathcal{V}_{\delta}(x)$, yields the following estimates:
\begin{align*}
\int_{0}^{1} \frac{1}{t^{d+1}} \int_{\mathcal{V}_{t\delta}(x)} |\nabla y (\xi)| |\xi - x|  d  \xi\, dt \leq  \lVert \nabla y \rVert_{L^{p}(\mathcal{V}_{\delta}(x))} \int_{0}^{1} \frac{1}{t^{d+1}}  \Big(\int_{\mathcal{V}_{t\delta}(x)} |\xi - x|^{q}  d   \xi\Big)^{\frac{1}{q}}dt \\
\leq  \lVert \nabla y \rVert_{L^{p}(\mathcal{V}_{\delta}(x))}  \int_{0}^{1} \frac{dt}{t^{\frac{d}{p}}} \Big(\int_{\mathcal{V}_{\delta}(x)} |x - z|^{q}  dz \Big)^{\frac{1}{q}} \leq \frac{\lVert \nabla y \rVert_{L^{p}(\Omega)}}{1-\frac{d}{p}}(\int_{\mathcal{V}_{\delta}(x)} |x - z|^{q}  dz \Big)^{\frac{1}{q}},
\end{align*}
where we used the boundedness of the $dt$-integral for $\frac{d}{p}<1$. Since $\rho(r): = \rho (|x-z|) = |x - z|^{q}$ is radially symmetric on $\mathcal{V}_{\delta}(x)$, integration in polar coordinates gives:
\begin{align*}
\Big(\int_{\mathcal{V}_{\delta}(x)} |x - z|^{q}  dz \Big)^{\frac{1}{q}} = |\mathcal{S}_{\delta}|^{\frac{1}{q}} \Big( \int_{0}^{\delta/2} \rho (r) r^{d-1} dr\Big)^{\frac{1}{q}} = C(d) \delta^{\frac{2d-1+q}{q}},
\end{align*}
where $\mathcal{S}_{\delta}$ is the $(d-1)$-sphere of radius $\frac{\delta}{2}$ with $|\mathcal{S}_{\delta}|= \frac{2 \pi^{d/2}}{\Gamma(d/2)}(\delta/2)^{d-1}$, $\Gamma$ is Euler's gamma function and $C(d) \in \mathbb{R}$ is $\delta$-independent. Dividing \eqref{MP: eq 7.2} by $|\mathcal{V}_{\delta}(x)|$ with $|\mathcal{V}_{\delta}(x)|=\frac{\pi^{d/2}}{\Gamma(d/2 + 1)}(\delta/2)^{d}$, combining the estimates and integrating, we get 
\begin{align*}
\left(\int_{\Omega_{\delta}}|y(x) - (Q_{\delta} y)(x) |^{2} \ dx \right)^{1/2} \leq  C_{p}(y) \delta^{\frac{d(2-q)+q-1}{q}}, 
\end{align*}
where $C_{p}(y): = C(d, p, \Omega) \lVert \nabla y \rVert_{L^{p}(\Omega)}$. The density of $C^{1}(\bar{\Omega})$ in $W^{1,p}(\Omega)$ extends the above estimate to the desired result.  
\end{proof}
\noindent
Under a mild assumption on the regularity of $y^{0} \in Y$, we prove the following result.
\begin{theorem}[\textbf{Upscaling consistency}] 
\label{Th: Upscaling consistency}
\textit{Suppose that Assumptions~\ref{Ass: Weak convergence}, ~\ref{Ass: Consistent parametrization} hold, $y^{0} \in Y \subset H^{2}(\Omega)$ with $\Omega \subset \mathbb{R}^{d}$, $d < p < \infty$ and $d\leq 3$. Then, for $\delta > 0$, it holds:
\begin{align*}
\underset{\varepsilon \rightarrow 0}{\lim} \ \mathcal{R}_{\delta}(y(u^{\varepsilon}), u^{\varepsilon}) \leq  C_{p}(y^{0})^{2} \delta^{\frac{2d(p-2)+2}{p}} +  2C_{\mathcal{L}[u]}^{2} + |\Omega \setminus \Omega_{\delta}|^{\frac{2+ \gamma}{\gamma}} C_{w},
\end{align*}
where $C_{p}(y^{0})\in \mathbb{R}_{+}$ and $C_{\mathcal{L}[u]}$ are according to Lemma~\ref{Lebesque approximation decay lemma} and Assumption~\ref{Ass: Consistent parametrization}, respectively, and $C_{w}\in \mathbb{R}_{+}$ is $\varepsilon$-independent.}
\end{theorem}
\begin{proof}  The triangle inequality and Young's inequality yield
\begin{align}\label{MP: eq 7.4}
\mathcal{R}_{\delta}(y(u^{\varepsilon}), u^{\varepsilon}) \leq 2\lVert \bar{Q}_{\delta} u^{\varepsilon} - y^{0} \rVert_{H}^{2} + 2 \lVert  y^{0} - y(u^{\varepsilon}) \rVert_{H}^{2}. 
\end{align}
Similarly, the estimation of the first term in \eqref{MP: eq 7.4} gives 
\begin{align}\label{MP: eq 7.5}
\lVert \bar{Q}_{\delta} u^{\varepsilon} - y^{0}  \rVert_{H}^{2} \leq 2\lVert Q_{\delta} u^{\varepsilon} - Q_{\delta}  y^{0} \rVert_{H^{\delta}}^{2} + 2 \lVert  y^{0} - Q_{\delta} y^{0} \rVert_{H^{\delta}}^{2} + \lVert u^{\varepsilon} -  y^{0} \rVert_{L^2(\Omega \setminus \Omega_{\delta})}^{2}. 
\end{align}
Clearly, Assumptions~\ref{Ass: Weak convergence} implies that $u^{\varepsilon}\rightharpoonup y^{0}$ as $\varepsilon \rightharpoonup 0$ in $H$. Thus, Lemma~\ref{L2 convergence lemma} guarantees that the first term on the right-hand side of \eqref{MP: eq 7.5} vanishes as $\varepsilon \to 0$. Since $\nabla y^{0} \in H^{1}(\Omega)$ due to $Y\subset H^2(\Omega)$ and $d\leq 3$,  the Sobolev embedding yields $\lVert \nabla y^{0} \rVert_{L^{p}(\Omega)} < \infty$ for $p\leq 6$. Observe further that $y^{0}$ is $\varepsilon$-independent, hence $C_{p}
(y^0)$ is also $\varepsilon$-independent. Therefore, Lemma~\ref{Lebesque approximation decay lemma} guarantees that $\lVert   y^{0} - Q_{\delta}  y^{0} \rVert_{H^{\delta}}^{2} \leq C_{p}(y^0)^{2} \delta^{\frac{2d(p-2)+2}{p}}$ and the latter estimate holds in the limit $\varepsilon \rightarrow 0$. 

We estimate the second term in \eqref{MP: eq 7.4} by appealing to Assumption~\ref{Ass: Consistent parametrization}. Next, let $s=\frac{2 + \gamma}{2}$ and $r=\frac{2+ \gamma}{\gamma}$. Then $\frac{1}{s} +\frac{1}{r}=1$ and the H\"older estimate holds
\begin{align*}
\lVert u^{\varepsilon} -  y^{0} \rVert_{L^2(\Omega \setminus \Omega_{\delta})}^{2} \leq |\Omega \setminus \Omega_{\delta}|^{\frac{2+ \gamma}{\gamma}} \lVert u^{\varepsilon} -  y^{0} \rVert_{L^{2+\gamma}(\Omega)}^{2} 
\end{align*}
for the third term in \eqref{MP: eq 7.5}. The uniform boundedness of $\{u^{\varepsilon}\}$ implies the existence of $C_{w}\in \mathbb{R}_{+}$ such that $\lVert u^{\varepsilon} -  y^{0} \rVert_{L^{2+\gamma}(\Omega)}^{2} \leq C_{w}$, hence completing the proof. 
\end{proof}
\noindent
Suppose that $\nabla y^{0} \notin L^{p}(\Omega)$ for some $p>d$. Then Lemma~\ref{Lebesque approximation decay lemma} is not applicable. In this case, the coupling term  \eqref{MP: eq 8} can be modified as follows: 
\begin{align}\label{Modified coupling term}
\mathcal{R}_{\delta}(y(u^{\varepsilon}), u^{\varepsilon}): = \lVert \bar{Q}_{\delta}u^{\varepsilon} - \bar{Q}_{\delta}y(u^{\varepsilon}) \rVert_{H}^{2}.    
\end{align}
Then, Theorem~\ref{Th: Upscaling consistency} still holds with $C_{p}(y^0)=0$.

\section{Learning-informed PDE-constrained optimization}\label{Section 3} We cast our hybrid physics-informed neural network based multiscale approach for a fixed $\varepsilon>0$ into the following learning-informed PDE-constrained optimization problem:
\begin{equation}\label{COCP: eq 1}
\left\{
\begin{split}
& \inf  J(y, v_{\boldsymbol{\theta},n}):=\mathcal{J}(v_{\boldsymbol{\theta},n}) + \tau_{2}\mathcal{R}_{\delta}(y,v_{\boldsymbol{\theta},n})  \ \text{over} \ (y, v_{\boldsymbol{\theta},n}) \in Y \times B_{\bar{r}^{\varepsilon},\theta, n}(u^{\varepsilon}),%\nonumber
\\ 
& \text{s.t.} \  e(y, v_{\boldsymbol{\theta},n})  = 0, %\ (y, v_{\boldsymbol{\theta},n}) \in Y \times \mathfrak{N}_{\theta, n}, %\nonumber
\end{split}\right.
\end{equation}
where $B_{\bar{r}^{\varepsilon},\theta, n}(u^{\varepsilon}):=B_{\bar{r}^{\varepsilon}}(u^\varepsilon) \cap \mathfrak{N}_{\theta, n}$ with $u^{\varepsilon}$, $\bar{r}^{\varepsilon}$ according to Assumption~\ref{Ass:Uniformity}, $J : Y\times U \rightarrow \mathbb{R}_{\geq 0}$ and fixed $\tau_{2}, \delta\in\mathbb{R}_{\geq 0}$. Further, let $e : Y\times U \rightarrow  Y^{\ast}$ be given by $e: (w,v)  \mapsto e(w,v) : =  b_{\mathcal{L}}[v](w, \cdot) - \langle f, \cdot \rangle_{Y^\ast,Y}$ and $\mathcal{R}_\delta$ denote the coupling term in \eqref{Modified coupling term}. Note that for $u^\varepsilon$ according to Assumption~\ref{Ass:Uniformity} there exists $v\in X$ arbitrarily close to $u^\varepsilon$ due to the density of the embedding $X\hookrightarrow U$. Then, Assumption~\ref{Ass:Uniform NN Approximation of elements} implies the existence of $\boldsymbol{\theta} \in \mathbb{R}^{N_{n_{\varepsilon}^{\ast}}}$ such that
\begin{align*}
\lVert u^{\varepsilon} - v_{\boldsymbol{\theta}, n_{\epsilon}^{\ast}} \rVert_{U} \leq \lVert u^{\varepsilon} -v  \rVert_{U} + C \lVert v - v_{\boldsymbol{\theta}, n_\varepsilon^\ast} \lVert_{X} \leq \bar{r}^{\varepsilon}, 
\end{align*}
for sufficiently large $n_{\varepsilon}^{\ast} \in \mathbb{N}$. Then, for all $n \geq n_{\varepsilon}^{\ast}$,  $B_{\bar{r}^{\varepsilon},\theta, n}(u^{\varepsilon}) \neq \emptyset$ and the following {\it fine-to-coarse scale map} is well-defined: 
\begin{align}\label{COCP: eq 2}
S : B_{\bar{r}^{\varepsilon},\theta, n}(u^{\varepsilon}) \subset U \rightarrow  Y, \quad u \mapsto y(u):=S(u),
\end{align}
with $e(y(u),u)=0$. Note that since \eqref{parameter iso} is non-convex in general, finding $\boldsymbol{\theta}$ with $v_{\boldsymbol{\theta},n} \in B_{\bar{r}^{\varepsilon},\theta, n}(u^{\varepsilon})$ requires solving an associated non-convex optimization problem. In practice, however, the latter can often not be guaranteed. 

We also need the following.
\begin{Assumption}[\textbf{Continuity}] 
\label{Ass: Continuity}
\textit{Let $u^{\varepsilon} \in U$ be the solution of \eqref{MP: eq 1} and $\{u_{k}^{\varepsilon}\}\subset B_{\bar{r}^{\varepsilon}}(u^{\varepsilon}) \cap X$ its approximating sequence. Then $S(u_{k}^{\varepsilon}) \rightharpoonup S(u^{\varepsilon})$ in $Y$ as $k \rightarrow \infty$.}
\end{Assumption}
\noindent
Eliminating $y$ from the set of independent variables in \eqref{COCP: eq 1} results in the reduced optimization problem 
\begin{equation} \label{NOC: eq 1}
%&\underset{v_{\boldsymbol{\theta},n} \in \mathfrak{N}_{\theta, n} \cap U }{
\inf \widehat{J}(v_{\boldsymbol{\theta},n}):= J(S(v_{\boldsymbol{\theta},n}), v_{\boldsymbol{\theta},n}) \ \    \text{over }v_{\boldsymbol{\theta},n} \in B_{\bar{r}^{\varepsilon},\theta, n}(u^{\varepsilon}).
\end{equation}
We note that guaranteeing the existence of minimizers in $\mathfrak{N}_{\theta, n}$ for \eqref{NOC: eq 1} is not possible, in general, as $\mathfrak{N}_{\theta, n}$ may not be topologically closed in $U$. %In fact it may not even be a subspace of $X$. 
In order to cope with this, we resort to the notion of quasi-minimization; cf. \cite{shin2023error}, see also \cite{brevis2022neural}. The latter only requires the existence of an infimum of $\widehat{J}$ over $B_{\bar{r}^{\varepsilon},\theta, n}(u^{\varepsilon})\ne\emptyset$. 

Clearly, since $\widehat{J}(\cdot) \geq 0$, $\inf \widehat{J}$ exists over $B_{\bar{r}^{\varepsilon},\theta, n}(u^{\varepsilon})$ for every $n \geq n_{\varepsilon}^{\ast}$. Now let $\{\gamma_n\}_{n \geq n_{\varepsilon}^{\ast}}$ be a real sequence with $\gamma_n>0$ for all $n \geq n_{\varepsilon}^{\ast}$ and $\gamma_n \downarrow 0$ as $n\to\infty$. Then, for every $n \geq n_{\varepsilon}^{\ast}$, there exists $u_{\boldsymbol{\theta},n}^{\varepsilon}\in B_{\bar{r}^{\varepsilon},\theta, n}(u^{\varepsilon})$ such that  
$$
\widehat{J}(u_{\boldsymbol{\theta},n}^{\varepsilon})\leq\inf_{v_{\boldsymbol{\theta},n}\in B_{\bar{r}^{\varepsilon},\theta, n}(u^{\varepsilon})}\widehat{J}(v_{\boldsymbol{\theta},n})+\gamma_n.
$$
We refer to $\{u_{\boldsymbol{\theta},n}^{\varepsilon}\}_{n \geq n_{\varepsilon}^{\ast}}$ as a sequence of quasi-minimizers of \eqref{NOC: eq 1}. In the following result, we consider $J$ with the coupling term \eqref{Modified coupling term}, but extending it to \eqref{MP: eq 8} is straightforward. 
\begin{theorem}\label{Loss convergence} Suppose that Assumptions~\ref{Ass:Uniformity}, ~\ref{Ass:Stability},  ~\ref{Ass:Uniform NN Approximation of elements}, and \ref{Ass: Continuity} hold. Let $\{u_{\boldsymbol{\boldsymbol{\theta}},n}^{\varepsilon}\}_{n \geq n_{\epsilon}^{\ast}}$, $ u_{\boldsymbol{\boldsymbol{\theta}},n}^{\varepsilon} \in  B_{\bar{r}^{\varepsilon},\theta, n}(u^{\varepsilon})$ be a quasi-minimizing sequence for $\widehat{J}:U\to\mathbb{R}_{\geq 0}$, where $B_{\bar{r}^{\varepsilon}}(u^{\varepsilon})$ is chosen according to Assumption~\ref{Ass:Uniformity} and $n_{\varepsilon}^{\ast} \in \mathbb{N}$ is chosen according to Assumption~\ref{Ass:Uniform NN Approximation of elements} to guarantee that $B_{\bar{r}^{\varepsilon},\theta, n} (u^{\varepsilon})\ne\emptyset$ for all $n\geq n_{\varepsilon}^{\ast}$. Then, $\underset{n \rightarrow \infty}{\lim}\widehat{J}(u_{\boldsymbol{\theta},n}^{\varepsilon}) \leq \tau_{2}\mathcal{R}_{\delta}(y(u^\varepsilon), u^{\varepsilon})$, where $u^\varepsilon$ is the solution to \eqref{MP: eq 1} according to Assumption~\ref{Ass:existence}. Moreover, for $\tau_{1} \geq 1$ it holds that
\begin{align}\label{a-posteriori bound}
\underset{n \rightarrow \infty}{\lim}\lVert u_{\boldsymbol{\theta}, n}^{\varepsilon} - u^{\varepsilon} \rVert_{V} \leq \frac{1}{\sqrt{C_{s}^{\varepsilon}}} \sqrt{\tau_{2} \mathcal{R}_{\delta}(y(u^{\varepsilon}), u^{\varepsilon})}.
\end{align}
\end{theorem}
\begin{proof}
Let $u^{\varepsilon}$ be the solution to \eqref{MP: eq 1} and $\{u_{k}^{\varepsilon}\}$, $u_{k}^{\varepsilon}\in X$ for all $k$, be its approximating sequence \eqref{MP: eq 3}. Let $\{r_{k}\}$ be a positive real sequence with $r_{k} \downarrow 0$ monotonically as $k \to\infty$ and choose $K_{\varepsilon}\in \mathbb{N}$ sufficiently large such that $B_{r_{k}}(u_{k}^{\varepsilon}):=\{v \in U : \ \lVert u_{k}^{\varepsilon} - v \rVert_{U} \leq  r_{k} \} \subset B_{\bar{r}^\varepsilon}(u^{\varepsilon})$ for all $k \geq K_{\varepsilon}$. Assumption~\ref{Ass:Uniform NN Approximation of elements} and the dense embedding $X \hookrightarrow U$ imply that $\lVert u_{k}^{\varepsilon} - v_{\boldsymbol{\theta}, n_{k}}^{\varepsilon}\rVert_{U}  \leq C \lVert u_{k}^{\varepsilon} - v_{\boldsymbol{\theta}, n_{k}}^{\varepsilon} \rVert_{X} \leq C\epsilon_{n_{k}} \leq r_{k} $ for $k\geq K_{\varepsilon}$ and some $v_{\boldsymbol{\theta}, n_{k}}^{\varepsilon} \in \mathfrak{N}_{\theta,  n_{k}}$, where $n_{k}:=n(k) \in \mathbb{N}$, $\epsilon_{n_{k}} \in \mathbb{R}_{+}$ and $C>0$ is some embedding constant. Observe further that $B_{r_{k}}(u_{k}^{\varepsilon}) \subset B_{\bar{r}^\varepsilon}(u^{\varepsilon})$ for $k\geq K_{\varepsilon}$ implies $\epsilon_{n_{K_{\varepsilon}}} \leq \epsilon_{n_{\varepsilon}^{\ast}}$. Therefore, $n_{K_{\varepsilon}} \geq n_{\varepsilon}^{\ast}$ with $\mathfrak{N}_{\theta,  n_{\varepsilon}^{\ast}} \subseteq \mathfrak{N}_{\theta,  n_{K_{\varepsilon}}}$ holds and there exists a sequence $\{v_{\boldsymbol{\theta}, n_{k}}^{\varepsilon}\}_{k\geq K_{\varepsilon}}$ with $v_{\boldsymbol{\theta}, n_{k}}^{\varepsilon} \in B_{r_{k}, \theta,  n_{k}}(u_{k}^{\varepsilon}) \subset B_{\bar{r}^\varepsilon}(u^{\varepsilon})$. We note that Assumption~\ref{Ass:Uniform NN Approximation of elements} implies that $\{\epsilon_{n_{k}}\}_{k \geq K_{\varepsilon}}$ converges to $0$ as $k\rightarrow \infty$, thus $n_{k} \rightarrow \infty$ as $k\rightarrow \infty$.  Then, once again the embedding $X \hookrightarrow U$ implies that
\begin{align*}
\lVert u^{\varepsilon} -  v_{\boldsymbol{\theta}, n_{k}}^{\varepsilon} \rVert_{U} \leq \lVert u^{\varepsilon} - u_{k}^{\varepsilon} \rVert_{U} + C \lVert u_{k}^{\varepsilon} - v_{\boldsymbol{\theta}, n_{k}}^{\varepsilon} \rVert_{X} \rightarrow 0 \quad \text{as} \quad k\rightarrow \infty.
\end{align*}
In addition, the boundedness of $\mathcal{A}^\varepsilon$ and $\mathcal{B}$, respectively, and the upper bound in \eqref{MP: subeq 2(a)} result in the existence of $C^\varepsilon\in \mathbb{R}_{+}$, possibly dependent on $\varepsilon$, such that
\begin{align*}
\lVert \mathcal{A}^{\varepsilon} v_{\boldsymbol{\theta}, n_{k}}^{\varepsilon} - f \rVert_{H} + \lVert \mathcal{B} v_{\boldsymbol{\theta}, n_{k}}^{\varepsilon}  \rVert_{Z} \leq C^\varepsilon \lVert v_{\boldsymbol{\theta}, n_{k}}^{\varepsilon} - u^{\varepsilon}_{k} \rVert_{X} +  \lVert \mathcal{A}^{\varepsilon} u_{k}^{\varepsilon} - f \rVert_{H} + \lVert \mathcal{B} u_{k}^{\varepsilon}  \rVert_{Z} \rightarrow 0
\end{align*}
as $k \rightarrow 0$. Therefore, $\{v_{\boldsymbol{\theta}, n_{k}}^{\varepsilon}\}_{k \geq K_{\varepsilon}}$ is also an approximating sequence \eqref{MP: eq 3}. 

Next, we study the coupling term and show that 
\begin{align}\label{Th 3.2. eq. 1}
\underset{k \rightarrow \infty}{\lim} \mathcal{R}_{\delta}(S(v_{\boldsymbol{\theta}, n_{k}}^{\varepsilon}), v_{\boldsymbol{\theta}, n_{k}}^{\varepsilon}) \leq \mathcal{R}_{\delta}(y(u^{\varepsilon}), u^{\varepsilon}).
\end{align}
Indeed, invoking the triangle inequality and $\bar{Q}_{\delta} \in L(H)$, we get 
\begin{align}\label{Th 3.2. eq. 2}
\mathcal{R}_{\delta}(S(v_{\boldsymbol{\theta}, n_{k}}^{\varepsilon}), v_{\boldsymbol{\theta}, n_{k}}^{\varepsilon}) \leq & \lVert \bar{Q}_{\delta} \rVert^{2} \lVert v_{\boldsymbol{\theta}, n_{k}}^{\varepsilon} - u_{k}^{\varepsilon} \rVert_{H}^{2} + \lVert \bar{Q}_{\delta}u_{k}^{\varepsilon} - \bar{Q}_{\delta} y(v_{\boldsymbol{\theta}, n_{k}}^{\varepsilon}) \rVert^{2}_{H} \\ 
+ & 2 \lVert \bar{Q}_{\delta} \rVert  \lVert v_{\boldsymbol{\theta}, n_{k}}^{\varepsilon} - u_{k}^{\varepsilon} \rVert_{H} \lVert \bar{Q}_{\delta}u_{k}^{\varepsilon} - \bar{Q}_{\delta} y(v_{\boldsymbol{\theta}, n_{k}}^{\varepsilon}) \rVert_{H} \nonumber. 
\end{align}
Since $v_{\boldsymbol{\theta}, n_{k}}^{\varepsilon} \in B_{\bar{r}^\varepsilon,\theta, n_{\varepsilon}^{\ast}}(u^{\varepsilon})$ for $k \geq K_{\varepsilon}$, the fine-to-coarse mapping \eqref{COCP: eq 2} is well-defined in \eqref{Th 3.2. eq. 2}. The embedding $U\hookrightarrow H$ guarantees that the first and the last term in \eqref{Th 3.2. eq. 2} vanish as $k\rightarrow 0$. Similarly, one estimates the intermediate term in \eqref{Th 3.2. eq. 2}: 
\begin{align}\label{Th 3.2. eq. 3}
\lVert \bar{Q}_{\delta} u_{k}^{\varepsilon} - \bar{Q}_{\delta} y (v_{\boldsymbol{\theta}, n_{k}}^{\varepsilon}) \rVert^{2}_{H}  & \leq \mathcal{R}_{\delta}(y(u^{\varepsilon}), u_{k}^{\varepsilon}) + \lVert \bar{Q}_{\delta}\rVert^{2} \lVert y(u^{\varepsilon}) -   y(v_{\boldsymbol{\theta}, n_{k}}^{\varepsilon})\rVert_{H}^{2} \\ 
&+ 2 \lVert \bar{Q}_{\delta} \rVert \lVert \bar{Q}_{\delta}  u_{k}^{\varepsilon} - \bar{Q}_{\delta} y(u^{\varepsilon}) \rVert_{H}\lVert y(u^{\varepsilon})-y(v_{\boldsymbol{\theta}, n_{k}}^{\varepsilon})\rVert_{H}.\nonumber
\end{align}
We use Assumption~\ref{Ass: Continuity} and the compactness of $Y \hookrightarrow H$ to obtain $\underset{k \rightarrow \infty}{\lim} \lVert y(v_{\boldsymbol{\theta}, n_{k}}^{\varepsilon}) - y(u^{\varepsilon}) \rVert_{H}^{2}=0$. Therefore, \eqref{Th 3.2. eq. 1} readily follows from \eqref{Th 3.2. eq. 2} and \eqref{Th 3.2. eq. 3}.

Next, we study the limit $\underset{n\rightarrow \infty}{\lim}\widehat{J}(u_{\boldsymbol{\theta}, n}^{\varepsilon})$. The stability bound \eqref{Ass:Stability} gives 
\begin{align}\label{Th 3.2. eq. 4}
\mathcal{J}(v_{\boldsymbol{\theta}, n_{k}}^{\varepsilon})\leq (\lVert f - \mathcal{A}^{\varepsilon} u_{k}^{\varepsilon} \rVert_{H} + C_{b}^\varepsilon r_{k})^{2}  + \tau_{1}(\lVert \mathcal{B} u_{k}^{\varepsilon} \rVert_{Z}  + C_{b}^\varepsilon r_{k})^{2}. 
\end{align} 
Since $B_{\bar{r}^\varepsilon, \theta,  n_{K_{\varepsilon}}}(u^{\varepsilon}) \neq \emptyset$  and $B_{\bar{r}^\varepsilon, \theta,  n_{K_{\varepsilon}}}(u^{\varepsilon}) \subset B_{\bar{r}^\varepsilon}(u^{\varepsilon})$, we can find a quasi-minimizer $u_{\boldsymbol{\theta}, n_{k}}^{\varepsilon} \in   B_{\bar{r}^\varepsilon, \theta, n_k}(u^{\varepsilon})$ for each $k \geq K_{\varepsilon}$ with $\gamma_{n_{k}}>0$ and $\gamma_{n_{k}} \downarrow 0$ such that
\begin{align} \label{Th 3.2. eq. 5}
\widehat{{J}}(u_{\boldsymbol{\theta}, n_{k}}^{\varepsilon}) \leq \underset{v \in B_{\bar{r}^\varepsilon, \theta, n_{k}}(u^{\varepsilon})}{\inf} \ \widehat{{J}}(v)  + \gamma_{n_k} \leq  \widehat{J}(v_{\boldsymbol{\theta}, n_{k}}^{\varepsilon})  + \gamma_{n_{k}}.
\end{align}
Then, it follows from \eqref{Th 3.2. eq. 1} and \eqref{Th 3.2. eq. 4} that $\underset{k \rightarrow \infty}{\lim} \widehat{J}(v_{\boldsymbol{\theta}, n_{k}}^{\varepsilon}) \leq  \tau_{2}\mathcal{R}_{\delta}(y(u^{\varepsilon}), u^{\varepsilon})$. Hence, for any $\epsilon_{J}>0$, there exists $K_{\epsilon_{J}} \geq K_{\varepsilon}$ such that $\widehat{J}(u_{\boldsymbol{\theta}, n_{k}}^{\varepsilon}) \leq \epsilon_{J}/2 + \tau_{2}\mathcal{R}_{\delta}(y(u^{\varepsilon}), u^{\varepsilon})$ for all $k \geq K_{\epsilon_{J}}$. By resorting to the notion of a quasi-minimizing sequence $\{u_{\boldsymbol{\boldsymbol{\theta}},n}^{\varepsilon}\}_{n \geq n_{\epsilon}^{\ast}}$, the existence of $N_{\epsilon_{J}} \geq n_{\varepsilon}^{\ast}$ follows with $\gamma_{n} \leq \epsilon_{J}/2$ for all  $n \geq N_{\epsilon_{J}}$. For $\hat{N}_{\epsilon_{J}} = \max \{n_{K_{\epsilon_{J}}}, N_{\epsilon_{J}}\}$ and $n \geq \hat{N}_{\epsilon_{J}}$, we have $\mathfrak{N}_{\theta, \hat{N}_{\epsilon_{J}}} \subseteq \mathfrak{N}_{\theta, n}$, and 
\begin{align*}
\widehat{J}(u_{\boldsymbol{\theta}, n}^{\varepsilon}) \leq \underset{v \in B_{\bar{r}^\varepsilon,\theta, n}(u^{\varepsilon})}{\inf} \widehat{J}(v) + \gamma_n \leq  \widehat{J}(u_{\boldsymbol{\theta}, \hat{N}_{\epsilon_{J}}}^{\varepsilon}) + \gamma_n \leq \epsilon_{J} + \tau_{2}\mathcal{R}_{\delta}(y(u^{\varepsilon}), u^{\varepsilon}).
\end{align*}
Since $\epsilon_{J}$ was arbitrarily chosen, the first limit claim is shown. 
 
The stability estimate \eqref{MP: subeq 2(a)} and $\tau_{1} \geq 1$ imply that 
\begin{align*}
C_{s}^{\varepsilon} \lVert u_{\boldsymbol{\theta}, n}^{\varepsilon} - u^{\varepsilon} \rVert_{V}^{2}  &\leq \lVert \mathcal{A}^{\varepsilon}u_{\boldsymbol{\theta}, n}^{\varepsilon} - \mathcal{A}^{\varepsilon}u^{\varepsilon} \rVert_{H}^{2} + \tau_{1}\lVert \mathcal{B} u_{\boldsymbol{\theta}, n}^{\varepsilon} - \mathcal{B} u^{\varepsilon} \rVert_{Z}^{2} \leq \widehat{J}(u_{\boldsymbol{\theta}, n}^{\varepsilon}).
\end{align*}
Then, \eqref{a-posteriori bound} follows from the previous result and the a-posteriori bound above. 
\end{proof}
\noindent
In our setting,  $\mathcal{R}_{\delta}$ serves the purpose of a regularizer in the training process, and one might expect from Theorem~\ref{Th: Upscaling consistency} and Theorem~\ref{Loss convergence} that its efficiency increases with decreasing $\varepsilon$. However, it may happen that $C_{s}^{\varepsilon} \rightarrow 0$ as $\varepsilon \rightarrow 0$ for the standard $H=L^{2}(\Omega)$ PINN objective, making the bound \eqref{a-posteriori bound} merely of qualitative nature. 
\begin{Remark}\label{PINN convergence results}
In general, $\inf \widehat{J}$ and $\inf \mathcal{J}$ over $B_{\bar{r}^\varepsilon,\theta, n}(u^{\varepsilon})$ are different. Therefore, quasi-minimizers of $\widehat{J}$ and $\mathcal{J}$ are also different. For a quasi-minimizing sequence $\{u_{\boldsymbol{\theta},n}^{\varepsilon}\}_{n \in \mathbb{N}}$ of $\mathcal{J}$, it holds that $ \lVert u_{\boldsymbol{\theta}, n}^{\varepsilon} - u^{\varepsilon} \rVert_{V} \leq \frac{1}{\sqrt{C_{s}^{\varepsilon}}} \sqrt{\mathcal{J}(u_{\boldsymbol{\theta},n}^{\varepsilon})}$ with $ \underset{n \rightarrow \infty}{\lim}\mathcal{J}(u_{\boldsymbol{\theta},n}^{\varepsilon}) = 0$; cf. \cite[Proposition 3.1 and Theorem 3.2]{shin2023error}.
\end{Remark}
\noindent

\subsection{Discrete approximation}

The (fully) discrete version of \eqref{COCP: eq 1} comes in three steps: (i) First, we consider $\boldsymbol{\theta}\to v_{\boldsymbol{\theta},n}$, reducing every NN-function to its finite set of generating parameters belonging to $\mathbb{R}^{N_{n}}$; (ii) then, we replace the state $y$ by a finite dimensional approximation $y_h$, with $h$ indicating the associated discretization parameter such as, e.g., the mesh width in a finite element method \cite{brenner2008mathematical}; (iii) finally, the discrete version $\hat{\mathfrak{J}}$ of $\hat{J}$ is replaced by its quadrature approximation $\hat{\mathfrak{J}}^{M, h}_{D}$, where $M$ indicates the number of collocation or quadrature points, and $D$ is related to a discretization of the compression operator \eqref{Compression operator}. 

Let us start with (i) and define
\begin{align}
\mathcal{S} :& = S \circ \mathfrak{F}_{n} :  \mathbb{R}^{N_{n}}  \rightarrow Y, \quad \boldsymbol{\theta}  \mapsto y(v_{\boldsymbol{\theta},n}):=\mathcal{S}(\boldsymbol{\theta}), \label{DOCP: eq 2 (a)}\\
\widehat{\mathfrak{J}}:& = \widehat{J} \circ \mathfrak{F}_{n}: \mathbb{R}^{N_{n}} \rightarrow \mathbb{R}_{\geq 0}, \quad \boldsymbol{\theta} \mapsto \widehat{\mathfrak{J}}(\boldsymbol{\theta}): = J(\mathcal{S}(\boldsymbol{\theta}), \mathfrak{F}_{n}(\boldsymbol{\theta})),\label{DOCP: eq 2 (b)}
\end{align}
where $\mathfrak{F}_{n}$ is according to \eqref{parameter iso}. 
We invoke the following assumption on differentiability.
\begin{Assumption}[\textbf{Fr\'echet differentiability}] 
\label{Ass: differentiability}
\textit{The fine-to-coarse scale map \eqref{COCP: eq 2} is continuously Fr\'echet differentiable and $\mathfrak{F}_{n} \in C^{\infty}(\mathbb{R}^{N_{n}}, X)$. }
\end{Assumption}
Assumption~\ref{Ass: differentiability} and the chain rule imply that $\mathcal{S}$ is continuously Fr\'echet differentiable. We also note that $\mathfrak{F}_{n}$ satisfies Assumption~\ref{Ass: differentiability} for smooth activation functions in a neural network, i.e., $\sigma \in C^{\infty}(\mathbb{R})$. Accordingly, the differentiability requirement can be reduced via reducing the one on $\sigma$. %For $v_{\boldsymbol{\theta},n} \in \mathfrak{N}_{\boldsymbol{\theta}, n}$ and with a slight mis-use of notation, we set $e(w, \boldsymbol{\theta}):=e(w, v_{\boldsymbol{\theta},n})$. 
Via the implicit function theorem for $e(y(v_{\boldsymbol{\theta},n}),\mathfrak{F}_n(\boldsymbol{\theta}))=0$ we obtain 
\begin{align} \label{DOCP: eq 4}
e_{y}(y(v_{\boldsymbol{\theta},n}), \boldsymbol{\theta})\mathcal{S}^{\prime}(\boldsymbol{\theta}) + e_{v}(y(v_{\boldsymbol{\theta},n}), \mathfrak{F}_n(\boldsymbol{\theta}))\mathfrak{F}_n'(\boldsymbol{\theta}) = 0,
\end{align}
where we have $\mathcal{S}^{\prime}(\boldsymbol{\theta})=S'(\mathfrak{F}_n(\boldsymbol{\theta}))\mathfrak{F}_n'(\boldsymbol{\theta}) \in L(\mathbb{R}^{N_{n}},Y)$ and, upon obvious identification, $\mathcal{S}^{\prime}(\boldsymbol{\theta})^*\in L(Y^*,\mathbb{R}^{N_{n}})$. 

By applying the chain rule, we find the gradient of \eqref{DOCP: eq 2 (b)}:
\begin{align}\label{DOCP: eq 5}
\nabla \widehat{\mathfrak{J}}(\boldsymbol{\theta})  = \mathcal{S}^{\prime}(\boldsymbol{\theta})^{*} \partial_{y}J (y(v_{\boldsymbol{\theta},n}),\mathfrak{F}_n(\boldsymbol{\theta})) +  \mathfrak{F}_n'(\boldsymbol{\theta})^{*}\partial_{v} J(y(v_{\boldsymbol{\theta},n}), \mathfrak{F}_n(\boldsymbol{\theta})),
\end{align}
In practical realizations of PINNs, the second summand in \eqref{DOCP: eq 5} is typically produced by automatic differentiation, and the first summand is realized via the adjoint method \cite{hinze2008optimization}. For the latter, we need the bilinear form $b_{\mathcal{L}^\ast}[v_{\boldsymbol{\theta}, n}](\cdot, \cdot): Y\times Y \rightarrow \mathbb{R}$ with
\begin{align}\label{adjoint_form}
b_{\mathcal{L}^\ast}[v_{\boldsymbol{\theta}, n}](w,v) : = \langle e_y (y(v_{\boldsymbol{\theta},n}),\mathfrak{F}_n(\boldsymbol{\theta}))^*w, v \rangle_{Y^\ast,Y}.
\end{align}
The following guarantees that \eqref{adjoint_form} is well-defined. 
\begin{Assumption}
\label{Ass: System invertibility}
\textit{Suppose that $B_{\bar{r}^\varepsilon, \theta, n_{\varepsilon}^{\ast}}(u^{\varepsilon}) \subset U$, where $\bar{r}^\varepsilon \in \mathbb{R}_{+}$ is according to Assumption~\ref{Ass:Uniformity} and $n_{\varepsilon}^{\ast}$ is chosen sufficiently large as in Assumption~\ref{Ass:Uniform NN Approximation of elements}. Then, there exist $C_{b^\ast}^{\varepsilon}, C_{c^\ast}^{\varepsilon} \in\mathbb{R}_+$ such that for all $v_{\boldsymbol{\theta},n} \in  B_{\bar{r}, \theta, n}(u^{\varepsilon})$ with $n \geq n_{\varepsilon}^{\ast}$ it holds that
\begin{align}
b_{\mathcal{L}^{\ast}}[v_{\boldsymbol{\theta},n}](v,w) \leq C_{b^\ast}^{\varepsilon} \lVert v \rVert_{Y} \lVert w \rVert_{Y}, \ \text{and} \ \  b_{\mathcal{L}^{\ast}}[v_{\boldsymbol{\theta},n}](v,v) \geq C_{c^\ast}^{\varepsilon} \lVert v \rVert^{2}_{Y}, \  \forall{v,w}\in Y.
\end{align}}
\end{Assumption}
%Assumption \ref{Ass: System invertibility} is equivalent to  $e_y(y(v_{\boldsymbol{\theta},n}),\boldsymbol{\theta})$ and its inverse being bounded for all $v_{\boldsymbol{\theta},n} \in B_{\bar{r}^\varepsilon, \theta, n}(u^{\varepsilon})$. 
Note that Assumption \ref{Ass: System invertibility} and \eqref{DOCP: eq 4} yield the adjoint equation
\begin{align*}
e_y(y(v_{\boldsymbol{\theta},n}),\mathfrak{F}_n(\boldsymbol{\theta}))^*p(v_{\boldsymbol{\theta},n})=-\partial_{y} J(y(v_{\boldsymbol{\theta},n}),\mathfrak{F}_n(\boldsymbol{\theta}))=2\tau_2\bar{Q}_{\delta}^*(\bar{Q}_{\delta}v_{\boldsymbol{\theta}, n}-\bar{Q}_{\delta}y(v_{\boldsymbol{\theta},n})),
\end{align*}
or its counterpart in weak form:
\begin{align}\label{DOCP: eq 6}
b_{\mathcal{L}^\ast}[v_{\boldsymbol{\theta}, n}](p(v_{\boldsymbol{\theta},n}),v)=2 \tau_{2} \langle \bar{Q}_{\delta}v_{\boldsymbol{\theta}, n} - \bar{Q}_{\delta}y(v_{\boldsymbol{\theta}, n}),\bar{Q}_{\delta}v \rangle_{H} \quad \forall{v} \in Y,
\end{align}
where $p(v_{\boldsymbol{\theta},n}):=p(y(v_{\boldsymbol{\theta},n}))\in Y$ denotes the adjoint variable (or adjoint state, sometimes also called co-state).

Now we come to the second step of discretization. Here we use the finite element (FE) method applied to the coarse-scale equation. More specifically, let $Y_{h}:=\text{span} \{ \phi_{j}, \ 1\leq j \leq N_{h}\} \subset Y$, $N_h\in\mathbb{N}$, be the standard finite dimensional space of piecewise-linear and globally continuous functions over a domain $\Omega\subset\mathbb{R}^d$. Of course, other choices of $Y_h$ are possible as well \cite{brenner2008mathematical}. The finite element approximation of \eqref{MP: eq 4}, which involves the neural network based function $v_{\boldsymbol{\theta}, n}$ as data, is then obtained by a standard Galerkin projection:
Find $y_{h}(v_{\boldsymbol{\theta}, n}) \in Y_{h}$ such that
\begin{align}\label{DOCP: state}
b_{\mathcal{L}}[v_{\boldsymbol{\theta}, n}](y_{h}(v_{\boldsymbol{\theta}, n}),v_{h}) = \langle \tilde{f}, v_{h} \rangle_{Y^{\ast}, Y} \quad \forall{v_{h}} \in Y_{h}. 
\end{align}
Assumption~\ref{Ass:Uniformity} and Assumption~\ref{Ass:Uniform NN Approximation of elements} imply that \eqref{DOCP: state} admits a unique solution $y_{h}(v_{\boldsymbol{\theta}, n}) \in Y_{h} $ for all $v_{\boldsymbol{\theta}, n} \in  B_{\bar{r}^\varepsilon, \theta, n}(u^{\varepsilon})$ with $n\geq n_{\varepsilon}^{\ast}$.  

The adjoint equation is discretized similarly: Find $p_{h}(v_{\boldsymbol{\theta},n}):=p_{h}(y_h(v_{\boldsymbol{\theta},n})) \in Y_{h}$ such that
\begin{align}\label{DOCP: adjoint}
b_{\mathcal{L}^\ast}[v_{\boldsymbol{\theta}, n}](p_{h}(v_{\boldsymbol{\theta},n}), v_h) = 2\tau_{2} \langle \bar{Q}_{\delta}v_{\boldsymbol{\theta}, n} - \bar{Q}_{\delta}y_{h}(v_{\boldsymbol{\theta},n}),\bar{Q}_{\delta}v_{h} \rangle_{H} \quad  \forall{v_{h}} \in Y_{h}.
\end{align}
Assumption~\ref{Ass: System invertibility} implies that \eqref{DOCP: adjoint} admits a unique solution $p_{h}(v_{\boldsymbol{\theta},n}) \in Y_{h}$  for $v_{\boldsymbol{\theta},n} \in B_{\bar{r}^\varepsilon, \theta, n}(u^{\varepsilon})$.  

Regarding both of the above equations he following error estimates hold true. 
\begin{theorem} Suppose that Assumptions~\ref{Ass:Stability}, ~\ref{Ass:Uniformity}, ~\ref{Ass:Uniform NN Approximation of elements}, ~\ref{Ass: Continuity} and ~\ref{Ass: System invertibility} hold. Furthermore, let $Y \subset H^{2}(\Omega)$ and for all $u_{1}, u_{2} \in B_{\bar{r}^\varepsilon}(u^\varepsilon)$ and $v,w \in Y$ it holds that
\begin{align}\label{bi_form_assumption}
|b_{\mathcal{L}}[u_{1}](v,w) - b_{\mathcal{L}}[u_{2}](v,w)| & \leq C^\varepsilon \lVert u_{1} - u_{2} \rVert_{V} \lVert v \rVert_{Y} \lVert w \rVert_{Y},
\end{align}
where $C^\varepsilon \in \mathbb{R}_{+}$. In addition, \eqref{bi_form_assumption} holds also for the adjoint form $b_{\mathcal{L}^{\ast}}[u](\cdot, \cdot)$. Then 
\begin{align}
\underset{h\rightarrow 0}{\lim} \  \underset{n\rightarrow \infty}{\lim}\lVert y(u^{\varepsilon}) - y_{h}(u_{\boldsymbol{\theta}, n}^{\varepsilon}) \rVert_{Y} &\leq C_{y}^{\varepsilon}\sqrt{\tau_{2}\mathcal{R}_{\delta}(y(u^{\varepsilon}), u^{\varepsilon})}, \label{state_covergence} \\
\underset{h \rightarrow 0}{\lim} \  \underset{n\rightarrow \infty}{\lim} \lVert p(u^{\varepsilon}) - p_{h}(u_{\boldsymbol{\theta}, n}^{\varepsilon}) \rVert_{Y} & \leq C_{ad}^{\varepsilon} \sqrt{\tau_{2}\mathcal{R}_{\delta}(y(u^{\varepsilon}), u^{\varepsilon})} \label{adjoint_covergence},
\end{align}
where $p(u^{\varepsilon}):=p(y(u^{\varepsilon})) \in Y$, $p_h(u_{\boldsymbol{\theta},n}^{\varepsilon}):=p_h(y_h(u_{\boldsymbol{\theta}, n}^{\varepsilon})) \in Y_h$, $C_{y}^{\varepsilon},C_{ad}^{\varepsilon} \in \mathbb{R}_{+}$, and the latter constants may depend on $\varepsilon$.
\end{theorem}
\begin{proof}
For $u_{\boldsymbol{\theta}, n}^{\varepsilon} \in B_{\bar{r}^\varepsilon, \theta, n}(u^{\varepsilon})$, we treat $b_{\mathcal{L}}[u_{\boldsymbol{\theta}, n}^{\varepsilon}](\cdot, \cdot)$ as an approximation of $b_{\mathcal{L}}[u^{\varepsilon}](\cdot, \cdot)$ and apply Strang's lemma \cite[Lemma 2.27]{alexandre2004theory} to get the estimate  
\begin{align}\label{Strang_estimate_state}
\lVert y(u^{\varepsilon}) - y_{h}(u_{\boldsymbol{\theta}, n}^{\varepsilon}) \rVert_{Y} \leq & e_{n,h}^{y, rhs} + \underset{v_{h} \in Y_{h}}{\inf}\Big( \big(1+\frac{C_{b}^{\varepsilon}}{C_{c}^{\varepsilon}}\big) \lVert y(u^{\varepsilon}) -v_{h}\rVert_{Y}  \\
& + \frac{1}{C_c^{\varepsilon}} \ \underset{w_h \in Y_h}{\sup}\frac{|b_{\mathcal{L}}[u^{\varepsilon}](v_h,w_h) - b_{\mathcal{L}}[u_{\boldsymbol{\theta}, n}^{\varepsilon}](v_h,w_h)|}{\lVert w_h\rVert_{Y}}\Big), \nonumber
\end{align}
where $e_{n,h}^{y, rhs}$ represents a discrepancy measure between the right-hand sides of \eqref{MP: eq 4} and \eqref{DOCP: state}, but $e_{n,h}^{y, rhs}=0$ in our case. Let $\mathcal{I}_{h}y(u^{\varepsilon})$ be the interpolant of  $y(u^{\varepsilon})$ in $Y_{h}$. Since $Y \subset H^{2}(\Omega)$, a well-known \cite{brenner2008mathematical} interpolation bound yields
\begin{align*}
\underset{v_{h} \in Y_{h}}{\inf} \lVert y(u^{\varepsilon}) - v_h \rVert_{Y} \leq \lVert y(u^{\varepsilon}) - \mathcal{I}_{h}y(u^{\varepsilon})\rVert_{Y} \leq Ch\lVert y(u^{\varepsilon}) \rVert_{H^{2}(\Omega)} \rightarrow 0 \quad \text{as} \quad h\rightarrow 0.
\end{align*}
Assumption~\ref{Ass:Uniformity} and the boundedness of $\mathcal{I}_{h}: Y \rightarrow Y_{h}$ imply that $\lVert \mathcal{I}_{h}y(u^{\varepsilon}) \rVert_{Y} \leq C$, where $C \in \mathbb{R}_{+}$ does not depend on $h$. Invoking \eqref{bi_form_assumption} and \eqref{a-posteriori bound}, we get
\begin{align*}
|b_{\mathcal{L}}[u^{\varepsilon}](\mathcal{I}_{h}y(u^{\varepsilon}),w_{h}) - b_{\mathcal{L}}[u_{\boldsymbol{\theta}, n}^{\varepsilon}](\mathcal{I}_{h}y(u^{\varepsilon}),w_{h})| & \leq \frac{C}{\sqrt{C_{s}^{\varepsilon}}}\sqrt{\tau_{2}\mathcal{R}_{\delta}(y(u^{\varepsilon}), u^{\varepsilon})} \lVert w_{h} \rVert_{Y}
\end{align*}
for $n \rightarrow \infty$. The result \eqref{state_covergence} then follows from \eqref{Strang_estimate_state} and the above estimates. 

The analysis of the adjoint equation is similar to that of the state equation. However, the discrepancy between the right-hand sides of \eqref{DOCP: eq 6} and \eqref{DOCP: adjoint}, i.e., the quantity
\begin{align*}
e_{n,h}^{p, rhs}: = \underset{v_h \in Y_h}{\sup} \frac{|2 \tau_{2} \langle \bar{Q}_{\delta} (u^{\varepsilon} - u_{\boldsymbol{\theta},n}^{\varepsilon} + y_{h}(u_{\boldsymbol{\theta}, n}^{\varepsilon})- y(u^{\varepsilon})), \bar{Q}_{\delta}v_h \rangle_{H}|}{\lVert v_h \rVert_{Y}}. 
\end{align*}
needs to be additionally estimated for the application of Strang's lemma. The Cauchy--Schwarz and triangle inequalities, the embeddings $Y \hookrightarrow H$, $V \hookrightarrow H$, and $\bar{Q}_{\delta} \in L(H)$ yield $C \in \mathbb{R}_{+}$ such that
\begin{align*}
e_{n,h}^{p, rhs} \leq C \bigg(\lVert u^{\varepsilon} - u_{\boldsymbol{\theta},n}^{\varepsilon} \rVert_V + \lVert y(u^{\varepsilon}) - y_h(u_{\boldsymbol{\theta}, n}^{\varepsilon}) \rVert_{Y}\bigg).
\end{align*}
Then, the result \eqref{adjoint_covergence} follows from \eqref{a-posteriori bound} and \eqref{state_covergence}. 
\end{proof}
Both FE discretized equations result in the following algebraic system:
\begin{align}\label{discrete_system_opt}
\mathbb{B}_{h}[\boldsymbol{\theta}] \boldsymbol{y}_{h} = \mathbb{F}_{h}, \quad \mathbb{B}_{h}[\boldsymbol{\theta}]^{\top} \boldsymbol{p}_{h} = 2\tau_{2}(\mathbb{P}_{h}[\boldsymbol{\theta}]- \mathbb{P}_{h}[\boldsymbol{y}_{h}]),
\end{align}
where $\boldsymbol{y}_{h} \in \mathbb{R}^{N_{h}}$ and $\boldsymbol{p}_{h} \in \mathbb{R}^{N_{h}}$ are the coefficients of the FE functions $y_{h}= \sum_{i=1}^{N_{h}}(\boldsymbol{y}_{h})_{i}\phi_{i}$ and $p_{h} = \sum_{i=1}^{N_{h}}(\boldsymbol{p}_{h})_{i} \phi_{i}$. Moreover, $\mathbb{B}_{h}[\boldsymbol{\theta}] \in \mathbb{R}^{N_{h} \times N_{h}}$,  $(\mathbb{B}_{h}[\boldsymbol{\theta}])_{ij}:=b_{\mathcal{L}}[v_{\boldsymbol{\theta}, n}](\phi_{i},\phi_{j})$, $\mathbb{F}_{h} \in \mathbb{R}^{N_{h}}$, $(\mathbb{F}_{h})_{j}:= \langle \tilde{f}, \phi_{j}\rangle_{Y^{\ast}, Y}$ and $\mathbb{P}_{h}[\boldsymbol{\theta}] \in \mathbb{R}^{N_{h}}$, $(\mathbb{P}_{h}[\boldsymbol{\theta}])_{j}:= \langle \bar{Q}_{\delta}v_{\boldsymbol{\theta}, n}, \bar{Q}_{\delta}\phi_{j} \rangle_{H}$.

For item (iii) we assume that $H=L^{2}(\Omega)$, $Z=L^{2}(\partial \Omega)$, and apply a Monte-Carlo approach to approximate the integrals in the PINN objective as follows:
\begin{align}\label{full_discrete_loss}
\mathcal{J}^{M}(v_{\boldsymbol{\theta},n}^{M}):= \frac{|\Omega|}{M_{\Omega}}\overset{M_{\Omega}}{\underset{i=1}{\sum}}\big(\mathcal{A}^{\varepsilon} v_{\boldsymbol{\theta},n}^{M}(x_{i}^{r}) - f^{\varepsilon}(x_{i}^{r})\big)^{2} +\frac{\tau_{1} |\partial \Omega|}{M_{\partial \Omega}}\overset{M_{\partial \Omega}}{\underset{i=1}{\sum}}\big(\mathcal{B}v_{\boldsymbol{\theta},n}^{M}(x_{i}^{b})\big)^{2},
\end{align}
where $\{x_{i}^{r}\}_{i=1}^{M_{\Omega}}$ and $\{x_{i}^{b}\}_{i=1}^{M_{\partial \Omega}}$ are (uniformly random) collocation points in $\Omega$ and   on $\partial \Omega$, respectively. Assuming that $M:=M_{\Omega} + M_{\partial \Omega}$ is sufficiently large to guarantee that $v_{\boldsymbol{\theta},n}^{M} \in B_{\bar{r}^\varepsilon, \theta, n}(u^{\varepsilon})$ for all $n \geq n_{\varepsilon}^{\ast}$, we discretize \eqref{Modified coupling term} using a quadrature rule: 
\begin{align}\label{full_discrete_couple_term}
\mathcal{R}_{\delta, D}^{h}(y_{h}(v_{\boldsymbol{\theta},n}^{M}), v_{\boldsymbol{\theta},n}^{M}): = \overset{N_{h}}{\underset{i=1}{\sum}} w_{i}^{h} \big((\bar{Q}_{\delta}^{D}v_{\boldsymbol{\theta}, n}^{M})(x_{i}^{h}) - (\bar{Q}_{\delta}^{D} y_{h}(v_{\boldsymbol{\theta},n}^{M}))(x_{i}^{h})\big)^{2},
\end{align}
where $\{x_{i}^{h}$, $w_{i}^{h}\}_{i=1}^{N_{h}}$ are our finite element nodes and quadrature weights, and 
\begin{align}\label{full_discrete_comp}
(\bar{Q}_{\delta}^{D}v)(x_{i}^{h}) = 
\begin{cases}
\frac{1}{|\mathcal{V}_{\delta}^{D}(x_{i}^{h})|} \underset{x_{j}^{h} \in \mathcal{V}_{\delta}^{D}(x_{i}^{h})}{\sum} v(x_{j}^{h}), \   x_{i}^{h} \in \Omega_{\delta}, \\
v(x_{i}^{h}),  \  x_{i}^{h} \in \Omega\setminus \Omega_{\delta},
\end{cases}
\end{align}
where $v \in C(\bar{\Omega})$ is a continuous function, $\mathcal{V}_{\delta}^{D}(x_{i}^{h}) := \{ x_{j}^{r} : \ |x_{j}^{r} - x_{i}^{h}| \leq  \frac{\delta}{2} \}$ corresponds to the set of residual collocation points around the coarse-scale mesh nodes, and $M_{D}$ is an associated number of such collocation points. We note that using \eqref{MP: eq 8} allows us not to apply $\bar{Q}_{\delta}^{D}$ to the coarse-scale approximation in \eqref{full_discrete_couple_term}, yielding a simple implementation. The discrete components \eqref{full_discrete_loss}, \eqref{full_discrete_couple_term} and \eqref{full_discrete_comp}, as well as the mapping \eqref{DOCP: eq 2 (a)}, finally yield the approximation $\widehat{\mathfrak{J}}^{M, h}_{D}$ of $\widehat{\mathfrak{J}}$. 

Algorithmically and assuming that $\nabla \widehat{\mathfrak{J}}^{M, h}_{D}(\boldsymbol{\theta})$ is a sufficiently accurate approximation of $\nabla \widehat{\mathfrak{J}}(\boldsymbol{\theta})$ (i.e. securing descent properties with respect to $\widehat{\mathfrak{J}}$ at $\boldsymbol{\theta}$) we use the discrete gradient, e.g., in the Adam optimizer \cite{kingma2014adam}, to minimize $\widehat{\mathfrak{J}}^{M, h}_{D}(\boldsymbol{\theta})$. Next we summarize our overall computational procedure in Algorithm~\ref{alg:hybrid learning}. Upon successful termination it produces neural network parameters and other outputs by setting $\hat{it}:= it + 1$.
\begin{algorithm}
  \caption{Hybrid physics-informed NN training}
  \label{alg:hybrid learning}
  \textbf{Input:} Tolerance $tol$, max. number of iterations $it_{\text{max}}$, initial NN parameters  $\boldsymbol{\theta}^{(0)}$, initial state FEM coefficients $\boldsymbol{y}_{h}(\boldsymbol{\theta}^{(0)})$, optimizer hyperparameters. \\
  \textbf{Output:} NN parameters $\hat{\boldsymbol{\theta}}:=\boldsymbol{\theta}^{(\hat{it})}$, control variable $v_{\hat{\boldsymbol{\theta}}, n}^{M} \approx u_{\boldsymbol{\theta}, n}^{\varepsilon, M}$, state variable $ \boldsymbol{y}_{h}(\hat{\boldsymbol{\theta}}) \in \mathbb{R}^{N_{h}}$. 
  \begin{algorithmic}[1]
   \WHILE { $0 \leq it \leq it_{\text{max}}-1$ or $\widehat{\mathcal{J}}^{M}(v_{\boldsymbol{\theta}^{(it)}, n}^{M}) > tol$}
  \STATE Solve the adjoint system $\mathbb{B}_{h}^{\textbf{ad}}[\boldsymbol{\theta}^{(it)}]  \boldsymbol{p}_{h}= 2\tau_{2}(\mathbb{P}_{h}[\boldsymbol{\theta}^{(it)}]- \mathbb{P}_{h}[ \boldsymbol{y}_{h}(\boldsymbol{\theta}^{(it)})])$
  %\STATE{Compute  $\boldsymbol{y}_{h}^{T} \mathbb{E}_{h}[\boldsymbol{\theta}_{k}^{(it)}] \ \boldsymbol{p}_{h}$  for $1\leq k \leq n$}
  \STATE{Compute $\nabla_{\theta} \widehat{J}^{M, h}_{D}(y_{h}(\boldsymbol{\theta}^{(it)}), \boldsymbol{\theta}^{(it)}):=\text{grad}_{\boldsymbol{\theta}}(\widehat{J}^{M, h}_{D}(\boldsymbol{\theta}^{(it)}))$}
  \STATE {Assemble the total gradient $\nabla \widehat{\mathfrak{J}}^{M, h}_{D}(\boldsymbol{\theta}^{(it)})$}
  \STATE {Update weights  $ \ \boldsymbol{\theta}^{(it+1)} \leftarrow \text{Optimizer}\big(\nabla \widehat{\mathfrak{J}}^{M, h}_{D}(\boldsymbol{\theta}^{(it)}), \ \text{optimizer hyperparameters} \big)$} 
  \STATE {Solve the state system $\mathbb{B}_{h}[\boldsymbol{\theta}^{(it+1)}] \boldsymbol{y}_{h}(\boldsymbol{\theta}^{(it+1)}) = \mathbb{F}_{h}$}
  \STATE {$it \leftarrow it+1$}
  \ENDWHILE
  \end{algorithmic}
\end{algorithm}
In Algorithm~\ref{alg:hybrid learning}, $\text{grad}_{\boldsymbol{\theta}}$ refers to automatic differentiation with respect to the NN parameters. The typical optimizer of choice also depends on hyperparameters. For the Adam optimizer, these include selecting a learning rate $lr \in \mathbb{R}_{+}$ or a schedule of learning rates $lr : \mathbb{N} \rightarrow \mathbb{R}_{+}$ with $it \mapsto lr(it)$, as well as specifying values for $\beta_{1}^{Ad} \in \mathbb{R}_{+}$ and $\beta_{2}^{Ad}\in \mathbb{R}_{+}$ for the moving average update, and setting a batch size. We initialize $\boldsymbol{\theta}^{(0)}$ using the Glorot scheme \cite{glorot2010understanding}, but $\boldsymbol{\theta}^{(0)}$ can also be obtained by solving a neighboring, yet simpler problem \cite{goswami2020transfer, xu2023transfer}. For an overview of standard NN optimization techniques, we refer exemplarily to \cite{goodfellow2016deep} and the references therein.

\section{An application in heat conduction}
\label{Section 4}
In this section, we focus on the example of heat conduction with respect to both, the coarse and the fine scale, respectively. For this purpose let $u_{\boldsymbol{\theta},n}$ denote the {\it parameterization} of $\mathcal{L}\big[u_{\boldsymbol{\theta},n}\big]$. Specifically, we study our hybrid approach in view of its embedded upscaling process.
%Afterwards, we apply the multiscale hybrid setup to the upscaling process.

\subsection{Upscaling-based parameterization}
Our exemplary stationary heat transfer problem in $\bar{\Omega}=[0,1]^{2}$ is defined as follows: 
\begin{align}\label{heat transfer problem}
-\nabla \cdot (\boldsymbol{K}^{\varepsilon} \nabla w^{\varepsilon})  = q,  \quad \text{in} \ \Omega, \quad\text{and}\quad w^{\epsilon}  = 0, \quad \text{on} \  \partial \Omega, 
\end{align}
where $w^{\varepsilon}$ is the temperature field, $q\in L^2(\Omega)$ is the source term, $\boldsymbol{K}^{\varepsilon}\in C^{0, 1}(\bar{\Omega},\mathbb{R}_+)$ is a Lipschitz continuous coefficient, and $\varepsilon \in \mathbb{R}_+$ is a small parameter indicating the fine scale length.  In addition, it holds that
\begin{align}\label{bounds on K}
\exists \alpha , \  \beta \in \mathbb{R}_{+} : \ \alpha \leq \boldsymbol{K}^{\varepsilon} (x) \leq \beta, \quad \forall{x} \in \bar{\Omega},
\end{align}
We assume further that $\lVert \nabla \boldsymbol{K}^{\varepsilon} \rVert_{L^{\infty}(\Omega)} \leq \frac{c_{\boldsymbol{K}}}{\varepsilon}$, where $c_{\boldsymbol{K}} \in \mathbb{R}_{+}$ is  $\varepsilon$-independent.

Solving $\eqref{heat transfer problem}$ using finite elements can be computationally demanding due to the requirement $h \ll \varepsilon$ for a mesh size $h\in \mathbb{R}_{+}$ and $\varepsilon \ll1$ in order to resolve the fine scale behaviour numerically. We aim to efficiently characterize material properties of $\Omega$ and determine a computationally feasible coarse-scale (homogenized) counterpart of $\eqref{heat transfer problem}$. In this vein, the concept of $G$-convergence is employed to formalize the notion of a homogenized equation and the related effective material; see, e.g., \cite{jikov2012homogenization}.
\begin{definition}
\label{Def: G-limit}
\textit{A coefficient sequence $\{\boldsymbol{K}^{\varepsilon}(\cdot)\}$ is said to $G$-converge to $\boldsymbol{K}^{*}(\cdot)$ as $\varepsilon \rightarrow 0$, if for any $q \in H^{-1}(\Omega)$ the sequence of solutions $\{w^{\varepsilon}\}$ of \eqref{heat transfer problem} satisfies}
\begin{align}\label{weak convergence in G limit}
w^{\varepsilon} \rightharpoonup w^{0} \quad \text{in} \quad H^{1}_{0}(\Omega), \quad \boldsymbol{K}^{\varepsilon}\nabla w^{\varepsilon} \rightharpoonup \boldsymbol{K}^{\ast}\nabla w^{0} \quad \text{in} \quad L^{2}(\Omega),
\end{align}   
\textit{where $w^{0}$ is the solution to the homogenized equation}
\begin{align*}
-\nabla \cdot (\boldsymbol{K}^{*} \nabla w^{0})  = q,  \quad \text{in} \ \Omega, \quad \text{and}\quad w^{0}  = 0, \quad \text{on} \  \partial \Omega. 
\end{align*}
\end{definition} 
\noindent
Various techniques exist to find the $G$-limit. Their respective applicability, however, depends on the specific problem and properties of the underlying medium. If such a $G$-limit exists, then it does not depend on $q$ and on the boundary data on $\partial \Omega$. We also note that the existence of $\boldsymbol{K}^{*}$ for general heterogeneous media still remains an open problem. Often, the representative volume element (RVE) technique can be applied to find an approximation of $\boldsymbol{K}^{*}$. This approach is widely utilized in engineering applications; see, e.g., \cite{christensen2012mechanics, durlofsky1991numerical, ewing2009simplified, farmer2002upscaling}. Here we state its equivalent characterization via weak convergence of gradients and fluxes in \eqref{weak convergence in G limit}: For any measurable set $\mathcal{V} \subseteq \Omega$ of measure $|\mathcal{V}|>0$ and $\langle \ \cdot \ \rangle_{\mathcal{V}} = \frac{1}{|\mathcal{V}|} \int_{\mathcal{V}} \cdot \ dx$ (understood component-wise), one has
\begin{align} \label{weak conv}
\underset{\varepsilon \rightarrow 0}{\lim} \  \langle \nabla w^{\varepsilon} \rangle_{\mathcal{V}} = \langle \nabla w^{0} \rangle_{\mathcal{V}}, \quad  \underset{\varepsilon \rightarrow 0}{\lim} \ \langle \boldsymbol{K}^{\varepsilon} \nabla w^{\varepsilon} \rangle_{\mathcal{V}} = \langle \boldsymbol{K}^{*} \nabla w^{0} \rangle_{\mathcal{V}}.
\end{align}
The existence of the above limits for general heterogeneous media is difficult to verify. Thus, here we only assume that these limits exist; see also Assumption~\ref{Ass: Weak convergence}. Further, we introduce the following approximations: 
\begin{align}
\langle \nabla w^{\varepsilon} \rangle_{\mathcal{V}} \approx  \underset{\varepsilon \rightarrow 0}{\lim} \ \langle \nabla w^{\varepsilon} \rangle_{\mathcal{V}}  = \langle \nabla w^{0} \rangle_{\mathcal{V}}, \ 
\langle \boldsymbol{K}^{\varepsilon} \nabla w^{\varepsilon} \rangle_{\mathcal{V}} \approx   \underset{\epsilon \rightarrow 0}{\lim} \ \langle \boldsymbol{K}^{\varepsilon} \nabla w^{\varepsilon} \rangle_{\mathcal{V}} \approx   \widetilde{\boldsymbol{K}} \langle \nabla w^{0} \rangle_{\mathcal{V}},  \nonumber
\end{align}
where $\widetilde{\boldsymbol{K}}$ is defined as follows. 
\begin{definition}[\textbf{Upscaled coefficient}]
\textit{The upscaled coefficient $\widetilde{\boldsymbol{K}}$ satisfies 
\begin{align}\label{constitutive_relation}
\langle \boldsymbol{K}^{\varepsilon} \nabla w^{\varepsilon} \rangle_{\mathcal{V}} = \widetilde{\boldsymbol{K}}(x) \langle \nabla w^{\varepsilon} \rangle_{\mathcal{V}},  \quad\text{for all }x\in \mathcal{V},
\end{align}
where $\widetilde{\boldsymbol{K}}$ is a $2\times 2$ tensor on $\mathcal{V}$ approximating $\boldsymbol{K}^{*}$ on $\mathcal{V}$.}
\end{definition}
Note that $\widetilde{\boldsymbol{K}}$ is constant on $\mathcal{V}$ as a consequence of its definition. Besides, it is beneficial to consider $\widetilde{\boldsymbol{K}}$ as a tensor for anisotropic materials. Further we observe that in our concrete setting, two solutions of the fine-scale problem are required to obtain $\widetilde{\boldsymbol{K}}$ from \eqref{constitutive_relation}. In contrast to $\boldsymbol{K}^{*}$, boundary conditions may then affect $\widetilde{\boldsymbol{K}}$. Of course, it is very desirable to identify a set of boundary conditions such that the dependence of $\widetilde{\boldsymbol{K}}$ on them is weak. In our application context, we choose the linear temperature drop boundary conditions $w_{i}^{\varepsilon}  = x_{i}$ on $\partial{\Omega}$, $i=1,2$, but periodic boundary conditions or temperature drop no-flow conditions can be applied as well. Their related analyses are similar; compare  \cite{durlofsky1991numerical, low2020influence, wu2002analysis}. The calculation of the upscaled thermal conductivity coefficient then leads to the following fine-scale problems: For $i \in \{1,2\}$, find $w_{i}^{\varepsilon}: \Omega \rightarrow \mathbb{R}$ with
\begin{align}\label{App.A: eq.1}
-\nabla \cdot (\boldsymbol{K}^{\varepsilon} \nabla w_{i}^{\varepsilon})  = q, \quad \text{in} \ \Omega, \quad\text{and}\quad w_{i}^{\varepsilon} = x_{i}, \quad   \text{on} \  \partial \Omega. 
\end{align}
We note that $q=0$ is typically chosen to prevent any influence of the source term on $\widetilde{\boldsymbol{K}}$. However, our multiscale solver is designed for computations of fine-scale solutions with $q \neq 0$, and setting $q=0$ defines a special case for our NN-based homogenization scheme below. Therefore, a more general form of \eqref{App.A: eq.1} is considered, but the influence of $q$ on $\widetilde{\boldsymbol{K}}$ is assumed to be rather weak. The upscaling process is described in Algorithm~\ref{alg:Upscaling algorithm}; cf. \cite[Section 4]{wu2002analysis}. 
\begin{algorithm}
  \caption{Upscaling algorithm}
  \label{alg:Upscaling algorithm}
  \textbf{Input:} $\Omega = \bigcup_{j=1}^{N} \mathcal{V}_{j}$ with $\mathcal{V}_{i} \cap \mathcal{V}_{j}=\emptyset$, $i \neq j$, the solutions $\boldsymbol{w}^{\varepsilon}=\{w_{1}^{\varepsilon}, w_{2}^{\varepsilon}\}$ of \eqref{App.A: eq.1}. \\
  \textbf{Output:} The upcaled coefficient $\boldsymbol{\widetilde{K}}$.
  \begin{algorithmic}[1]
  \FOR{$j \in \{1,\dots,N\}$}
  \STATE Compute $\boldsymbol{\widetilde{F}_{i}^{j}}: =  \langle \boldsymbol{K}^{\varepsilon} \nabla w_{i}^{\varepsilon} \rangle_{\mathcal{V}_{j}}\in \mathbb{R}^{2}, \ \boldsymbol{\widetilde{T}_{i}^{j}}: = \langle \nabla w_{i}^{\varepsilon} \rangle_{\mathcal{V}_{j}} \in \mathbb{R}^{2}$ for $i\in \{1,2\}$.
  \STATE Insert $\boldsymbol{\widetilde{F}_{i}^{j}}$ and $\boldsymbol{\widetilde{T}_{i}^{j}}$ into \eqref{constitutive_relation} and find $\boldsymbol{\widetilde{K}}^{j}\in \mathbb{R}^{2 \times 2}$  from the matrix equation:
  \small{
    \begin{align}\label{App.A: eq.2}
    \begin{pmatrix}
        \boldsymbol{\widetilde{K}}^{j}_{11} \ \boldsymbol{\widetilde{K}}^{j}_{12} \\
        \boldsymbol{\widetilde{K}}^{j}_{21} \ \boldsymbol{\widetilde{K}}^{j}_{22}
    \end{pmatrix} 
    \begin{pmatrix}
        (\boldsymbol{\widetilde{T}_{1}^{j}})_{1} \ ( \boldsymbol{\widetilde{T}_{2}^{j}})_{1} \\
        (\boldsymbol{\widetilde{T}_{1}^{j}})_{2} \ ( \boldsymbol{\widetilde{T}_{2}^{j}})_{2}
    \end{pmatrix} = \begin{pmatrix}
       (\boldsymbol{\widetilde{F}_{1}^{j}})_{1}  \  (\boldsymbol{\widetilde{F}_{2}^{j}})_{1} \\
       (\boldsymbol{\widetilde{F}_{1}^{j}})_{2}  \  (\boldsymbol{\widetilde{F}_{2}^{j}})_{2}
    \end{pmatrix}.
    \end{align}}
    \ENDFOR
  \STATE Get $\boldsymbol{\widetilde{K}}$ with $\boldsymbol{\widetilde{K}}(x)=\boldsymbol{\widetilde{K}}^{j}$ for $ x \in \mathcal{V}_{j}$
  \end{algorithmic}
\end{algorithm}

Clearly, we need to ensure that the matrix of averaged gradients in \eqref{App.A: eq.2} is invertible. This requires the partition $\{\mathcal{V}_{j}\}$ to be sufficiently heterogeneous with each $\mathcal{V}_{j}$ of ``reasonable'' size to prevent linear dependence between $\boldsymbol{\widetilde{T}_{1}^{j}}$ and $\boldsymbol{\widetilde{T}_{2}^{j}}$. Algorithm~\ref{alg:Upscaling algorithm} is based on problem-dependent constitutive relations. In our case, we apply Fourier's law of heat conduction \eqref{App.A: eq.2}; cf. \cite{ewing2009simplified, low2020influence}. In view of other applications, Darcy's law is used for porous media flows \cite{griebel2010homogenization, iliev2013numerical, wu2002analysis} as the relation between fluid velocity and pressure, and Hooke's law is applied \cite{chalon2004upscaling} in elasticity to relate stress and strain fields.   

For $\mathcal{V}=\Omega=(0,1)^{2}$, we derive a simplified formula for the upscaled coefficient; cf. \cite{wu2002analysis}. For this purpose, it is convenient to study problem \eqref{App.A: eq.1} as a problem with homogeneous Dirichlet boundary conditions. Let $u^{x_{i}}:=x_{i}$ be the extension of our boundary conditions to $\Omega$ and consider the problem
\begin{align}\label{App.A: eq.3} 
-\nabla \cdot (\boldsymbol{K}^{\varepsilon} \nabla u_{i}^{\varepsilon})  = q  +   \nabla \cdot (\boldsymbol{K}^{\varepsilon} \nabla u^{x_{i}}), \quad \text{in} \ \Omega, \quad\text{and}\quad u_{i}^{\varepsilon} = 0, \quad   \text{on} \  \partial \Omega. 
\end{align}
Now, we apply Algorithm~\ref{alg:Upscaling algorithm} to $\boldsymbol{\widetilde{F}_{i}}: =  \langle \boldsymbol{K}^{\varepsilon} \nabla w_{i}^{\varepsilon} \rangle_{\Omega}$ and $\boldsymbol{\widetilde{T}_{i}}: =  \langle   \nabla w_{i}^{\varepsilon} \rangle_{\Omega}$. Since $w_{i}^{\varepsilon} = u_{i}^{\varepsilon} + u^{x_{i}}$ and $\nabla u^{x_{i}} = \boldsymbol{e}_{i}$, where $\boldsymbol{e}_{i}\in\mathbb{R}^2$ denotes the $i$-th unit vector, the divergence theorem yields
\begin{align}\label{App.A: eq.4}
\langle \nabla w_{i}^{\varepsilon} \rangle_{\Omega} = \int_{\partial \Omega} u_{i}^{\varepsilon} \ \boldsymbol{\eta} \ ds + \int_{\Omega} \boldsymbol{e}_{i} \  dx  = \boldsymbol{e}_{i},
\end{align}
where $\boldsymbol{\eta}(x)$ is the outward unit normal at $x \in \partial \Omega$. Inserting \eqref{App.A: eq.4} into \eqref{App.A: eq.2}, we get 
\begin{align}\label{App.A: eq.5}
\boldsymbol{\widetilde{K}}_{ij} =  \int_{\Omega}\boldsymbol{K}^{\varepsilon} \partial_{x_{j}}  w_{i}^{\varepsilon} \ dx = \int_{\Omega}\boldsymbol{K}^{\varepsilon} (\partial_{x_{j}} u_{i}^{\varepsilon}  + \delta_{ij}) \ dx,
\end{align}
where $\delta_{ij}$ denotes the Kronecker delta symbol and $\partial_{x_{j}}:= \frac{\partial}{\partial_{x_{j}}}$.

\subsection{Analysis of the model equations}\label{Analysis of the model equations} From now on we assume that the following assumption is satisfied. 
\begin{Assumption}
\label{Ass: Coefficient form}
Assume that $u^{\varepsilon}: = u^{\varepsilon}_{1}$ is the solution to problem \eqref{App.A: eq.3}. Moreover, $\boldsymbol{\widetilde{K}}[u^{\varepsilon}]:= \boldsymbol{\widetilde{K}} \in \mathbb{R}^{2\times 2}$, where $\boldsymbol{\widetilde{K}}_{12}=\boldsymbol{\widetilde{K}}_{21}\geq 0$ and $\boldsymbol{\widetilde{K}}_{22}\geq 0$ are provided to us and fixed, and $\boldsymbol{\widetilde{K}}[u^{\varepsilon}]_{11}= \int_{\Omega} \boldsymbol{K}^{\varepsilon}(\partial_{x_{1}} u^{\varepsilon} + 1 ) \ dx $ is according to \eqref{App.A: eq.5}. 
\end{Assumption}
Assumption \ref{Ass: Coefficient form} implies that one needs to solve the model fine-scale problem \eqref{App.A: eq.3} only for $i=1$ to get $\boldsymbol{\widetilde{K}}_{11}$. We note that for periodic coefficients $\boldsymbol{K}^{\varepsilon}$ with period $\varepsilon$, $\boldsymbol{\widetilde{K}}_{12}=\boldsymbol{\widetilde{K}}_{21}=0$  and $\boldsymbol{\widetilde{K}}[u^{\varepsilon}] = \boldsymbol{\widetilde{K}}[u^{\varepsilon}]_{11} I$, where $I \in \mathbb{R}^{2\times 2}$ is the identity. 

Next, we show that the assumptions of Section~\ref{Section 2} are satisfied by our fine- and coarse-scale model problems. For this purpose, let $X:= C^{k}(\bar{\Omega})$ for some $k\geq 2$, $U:=H^{2}(\Omega)$, $U_0:=H^{2}(\Omega) \cap H_{0}^{1}(\Omega)$ and $V:=H^{1/2}(\Omega)$ be equipped with the norms
\begin{align*}
\lVert v\rVert_{X}= & \underset{|\mathfrak{a}|\leq 2} {\sum}  \ \underset{x \in \Omega}{\sup} \ |\partial_{x}^{\mathfrak{a}} v(x)|, \quad \lVert v \rVert_{U}=\lVert v \rVert_{U_0} = \big(\underset{|\mathfrak{a}|\leq 2}{\sum} \lVert \partial_{x}^{\mathfrak{a}} v \rVert_{L^{2}(\Omega)}^{2} \big)^{1/2}, \\
& \lVert v \rVert_{V} = \left(\lVert v \rVert_{L^{2}(\Omega)}^{2} + \int_{\Omega} \int_{\Omega} \frac{|v(x) - v(y)|^{2}}{|x-y|^{d+1}} \ dx dy \right)^{1/2}, 
\end{align*}
respectively. Recall that $X$ is dense in $U$ and $\lVert u\rVert_{U} \leq |\Omega|^{1/2} \lVert u \rVert_{X}$ for $v \in X$, i.e, $X \hookrightarrow U$. Let $H:=L^{2}(\Omega)$, $Z:=L^{2}(\partial \Omega)$, $Y:=H^{1}_{0}(\Omega)$, $Y^{*} = H^{-1}(\Omega)$ with $Y \hookrightarrow H \hookrightarrow Y^{*}$, where $Y$ is compactly embedded into $H$ by the Rellich--Kondrachov Theorem. We define the following linear operators 
\begin{align*}
\mathcal{A}^{\varepsilon}: U \rightarrow H,  \quad   v \mapsto \mathcal{A}^{\varepsilon}v=- \nabla \cdot (\boldsymbol{K}^{\varepsilon} \nabla v), \ \text{and} \ \mathcal{B}: U \rightarrow Z,  \quad v\mapsto \mathcal{B}v = \restr{v}{\partial \Omega}.
\end{align*}
%For the discussion below, the sub-index $i \in \{1, 2\}$ is fixed for the fine-scale solution $u_{i}^{\varepsilon}$ and the coarse-scale solution $y_{i}$.
We cast our model fine-scale problem into the operator form
\begin{align} \label{HC eq 1}
\mathcal{A}^{\varepsilon} u  = f^{\varepsilon},   \quad \text{in} \ H, \quad \text{and}\quad \mathcal{B}u = 0, \quad \text{in} \  Z, 
\end{align}
with $f^{\varepsilon}: = q  -   \mathcal{A}^{\varepsilon}  u^{x_{1}} \in L^2(\Omega)=H$. By the Lax--Milgram lemma, \eqref{HC eq 1} has a unique solution $u^{\varepsilon} \in Y$, and we have $w_{1}^{\varepsilon}=u^{\varepsilon} + u^{x_{1}}$. The next proposition shows that $u^{\varepsilon} \in U$, thereby confirming Assumption~\ref{Ass:existence}, and that  Assumption~\ref{Ass:Stability} is satisfied.

\begin{proposition}\label{Stability independence} There exist $C_{s}^{\varepsilon} \in \mathbb{R}_{+}$ and $C_{b}^{\varepsilon} \in \mathbb{R}_{+}$ such that
\begin{align}\label{stability_bound 2}
C_{s}^{\varepsilon} \lVert u \rVert_{U}^{2} \leq  \lVert \mathcal{A}^{\varepsilon}u \rVert_{H}^{2} \leq C_{b}^{\varepsilon} \lVert u \rVert_{U}^{2},  \quad \forall{u} \in U_{0},
\end{align}
where $C_{s}^{\varepsilon} \rightarrow 0$, $C_{b}^{\varepsilon} \rightarrow \infty$ as $\varepsilon \rightarrow 0$.
\end{proposition}
\begin{proof}
We multiply \eqref{HC eq 1} (with $u=u^{\varepsilon}$) by $u^{\varepsilon}$ and integrate by parts to get
\begin{align*}  
\alpha \int_{\Omega}|\nabla u^{\varepsilon} |^{2} \ dx \leq \int_{\Omega} \boldsymbol{K}^{\varepsilon}  \nabla u^{\varepsilon} \cdot \nabla u^{\varepsilon} \ dx  =  \int_{\Omega}   q u^{\varepsilon} - \boldsymbol{K}^{\varepsilon} \partial_{x_{1}} u^{\varepsilon} \ dx. 
\end{align*}
Using the Cauchy-Schwarz and Poincaré inequalities, we get
\begin{align}\label{grad_bound}
\lVert \nabla u^{\varepsilon} \rVert_{H} \leq \frac{\beta + c_{p} \lVert q \rVert_{H}}{\alpha},
\end{align}
where the Poincaré constant $c_{p}\in \mathbb{R}_{+}$ depends only on $\Omega$, and $\beta$ is according to \eqref{bounds on K}. Note that problem \eqref{HC eq 1} can be written as follows:
\begin{align*}
- \Delta u^{\varepsilon} = g^{\varepsilon}:=\frac{f^{\varepsilon} + \nabla \boldsymbol{K}^{\varepsilon} \cdot \nabla u^{\varepsilon}}{\boldsymbol{K}^{\varepsilon}},  \quad \text{in} \ \Omega, \quad \text{and}\quad  u^{\varepsilon} = 0, \quad    \text{on} \  \partial \Omega \nonumber.
\end{align*}
Invoking a standard $H^{2}(\Omega)$ regularity result for convex domains \cite{grisvard2011elliptic}, we get 
\begin{align}\label{H2-regularity}
\lVert u^{\varepsilon} \rVert_{U} \leq \widehat{C} \lVert g^{\varepsilon} \rVert_{H} \leq \frac{\widehat{C}}{\alpha} \big(\lVert q\rVert_{H} + (\lVert \nabla u^{\varepsilon} \rVert_{H}+1)\lVert  \nabla \boldsymbol{K}^{\varepsilon} \rVert_{L^{\infty}(\Omega)}\big).
\end{align} 
where $\widehat{C}=\widehat{C}(\Omega) \in \mathbb{R}_{+}$. Using \eqref{grad_bound} and $\lVert \nabla \boldsymbol{K}^{\varepsilon} \rVert_{L^{\infty}(\Omega)} \leq \frac{c_{\boldsymbol{K}}}{\varepsilon}$, we obtain the estimate 
\begin{align}\label{bounded_inverse}
\lVert u^{\varepsilon} \rVert_{U} \leq \frac{\widehat{C}}{\alpha} \big((1+\frac{\beta}{\alpha}) \frac{c_{\boldsymbol{K}}}{\varepsilon}  + (1 + \frac{c_{\boldsymbol{K}}}{\varepsilon}\frac{c_{p}}{\alpha}) \lVert q \lVert_{H}\big).
\end{align}
The bound \eqref{bounded_inverse} implies that $\lVert (\mathcal{A}^{\varepsilon})^{-1} \rVert \leq (C_{1} + \frac{C_{2}}{\varepsilon}) < \infty$ for $\varepsilon > 0$ and some $\varepsilon$-independent $C_1, C_{2} \in \mathbb{R}_{+}$, from where the lower bound in \eqref{stability_bound 2} readily follows. Using the Young and the Cauchy--Schwarz inequalities, for $u \in U$, we get
\begin{align}\label{upper bound hc}
\lVert \mathcal{A}^{\varepsilon} u \rVert_{H}^{2} & \leq 2 \big( \lVert \nabla \boldsymbol{K}^{\varepsilon} \cdot \nabla u \rVert_{H}^{2}  + \lVert \boldsymbol{K}^{\varepsilon}  \Delta u \rVert_{H}^{2}\big) \leq 
\left(\frac{C}{\varepsilon^2}+\beta^2\right) \lVert u \rVert_{U}^{2},
\end{align}
where $C\in \mathbb{R}_{+}$ is $\varepsilon$-independent. This proves the upper bound in \eqref{stability_bound 2}.
\end{proof}
\noindent
The estimate \eqref{stability_bound 2} and interpolation techniques are used to verify the lower bound in Assumption~\ref{Ass:Stability} including $ \lVert \mathcal{B}u \rVert_{Z}$, see \cite[Lemma 4.1]{zeinhofer2023unified} and references therein. The upper bound in Assumption~\ref{Ass:Stability} follows from \eqref{upper bound hc} and the trace inequality. These estimates are summarized in the following proposition. 
\begin{proposition}\label{Stability soft penalty} There exist $\bar{C}_{s}^{\varepsilon} \in \mathbb{R}_{+}$ and $C_{b}^{\varepsilon} \in \mathbb{R}_{+}$ such that
\begin{align}\label{stability_bound 3}
\bar{C}_{s}^{\varepsilon} \lVert u \rVert_{V}^{2} \leq \lVert \mathcal{A}^{\varepsilon}u \rVert_{H}^{2} + \lVert \mathcal{B} u \rVert_{Z}^{2} \leq C_{b}^{\varepsilon} \lVert u \rVert_{U}^{2}, \quad \forall{u} \in U, 
\end{align}
where $\bar{C}_{s}^{\varepsilon} \rightarrow 0$, $C_{b}^{\varepsilon} \rightarrow \infty$ as $\varepsilon \rightarrow 0$.
\end{proposition}

The next result \cite[Theorem 3.3]{xie2011errors} shows that NNs with smooth activation functions are universal approximators of $C^{k}(\bar{\Omega})$, verifying Assumption~\ref{Ass:Uniform NN Approximation of elements} for $X:=C^{2}(\bar{\Omega})$. 
\begin{proposition}\label{UA smooth} Suppose that $\sigma \in C^{\infty}(\mathbb{R}), \  \sigma^{(s)}(0) \neq 0$ for $s=0,1,...,$ and $\bar{\Omega}:=[0,1]^{d}$. If $v \in C^{k}(\bar{\Omega})$, then there exists an architecture $\vec{\boldsymbol{n}}$ with one hidden layer and  $v_{\boldsymbol{\theta}, n} \in \mathfrak{N}_{\theta, n}$ such that
\begin{align*}
\underset{x \in \Omega}{\sup} \ |\partial_{x}^{\mathfrak{a}}v(x) - \partial_{x}^{\mathfrak{a}} v_{\boldsymbol{\theta}, n}(x)|  = \mathcal{O}\left( \frac{1}{n^{(k-|\mathfrak{a}|)/d}} \omega \big(\partial_{x}^{\mathfrak{b}}v, \frac{1}{n^{1/2}} \big) \right)
\end{align*}
holds for all multi-indices $\mathfrak{a}, \mathfrak{b}$ with $|\mathfrak{a}|\leq k, |\mathfrak{b}|=k$, where, for $\delta_c>0$, $\omega(v, \delta_{c}) = \sup \{|v(x)-v(y)|:|x-y| \leq \delta_{c}, x, y \in \bar{\Omega}\}$ is the modulus of continuity of $v$.
\end{proposition}
\noindent
For example, the swish function $\sigma(x)=x \ \text{sigmoid}(x)$ satisfies the prerequisites of Proposition~\ref{UA smooth}, but the widely-used $\tanh(x)$ does not satisfy $\sigma^{(2)}(0) \neq 0$ rendering Proposition~\ref{UA smooth} non-applicable. The following result \cite[Theorem 5.1]{de2021approximation} is then useful.
\begin{proposition}\label{tanh approximation} Suppose that $\sigma(x)=\tanh(x)$ and $\bar{\Omega}:=[0,1]^{d}$. If $v \in C^{k}(\bar{\Omega})$, then there exists an architecture $\vec{\boldsymbol{n}}$ with two hidden layers and $v_{\boldsymbol{\theta}, n} \in \mathfrak{N}_{\theta, n}$ such that $\lVert v - v_{\boldsymbol{\theta}, n} \rVert_{W^{m, \infty}(\bar{\Omega})}= \mathcal{O}(n^{-(k-m)})$ holds\footnote{We state the asymptotic convergence rate, but explicit constants are estimated in \cite{de2021approximation}.}, where $\lVert w\rVert_{W^{m, \infty}(\bar{\Omega})}= \underset{|\mathfrak{a}|\leq m}{\max} \ \underset{x \in \Omega}{\sup} \ |\partial_{x}^{\mathfrak{a}} w(x)|$.
\end{proposition}
\noindent
From the prerequisites of Proposition~\ref{tanh approximation} we infer $X:=C^{k}(\bar{\Omega})$ with $k>2$ as $m=2$ in our example. Clearly, $\mathfrak{N}_{\theta, n} \subset X$ and $\lVert v \rVert_{X}\leq C \lVert v \rVert_{W^{2, \infty}(\bar{\Omega})}$ for $v \in X$ and some $C \in \mathbb{R}_{+}$. The latter implies $X\subset\overline{\cup_{n}\mathfrak{N}_{\theta, n}}$, thereby verifying Assumption~\ref{Ass:Uniform NN Approximation of elements} for the hyperbolic tangent activation function. 

Our model coarse-scale problem in its weak form reads: Find $y(u^{\varepsilon}) \in Y$ such that
\begin{align}\label{weak formulation}
b_{\mathcal{L}}[u^{\varepsilon}](y(u^{\varepsilon}),v) = \langle \tilde{f}, v \rangle_{Y^{\ast}, Y} \quad \forall{v} \in Y,
\end{align}
where $\tilde{f}: = q + \nabla \cdot (\boldsymbol{\widetilde{K}}[u^{\varepsilon}] \nabla u^{x_{1}})$ and the forms are defined as follows:
\begin{align*}
 b_{\mathcal{L}}[u^{\varepsilon}](w,v): = \int_{\Omega} \boldsymbol{\widetilde{K}}[u^{\varepsilon}] \nabla w \cdot \nabla v \ dx, \
\langle \tilde{f}, v \rangle_{Y^{\ast},Y}  = \int_{\Omega} q v \ dx -  \int_{\Omega} \widetilde{\boldsymbol{K}}[u^{\varepsilon}] \boldsymbol{e}_{1} \cdot \nabla v \ dx.
\end{align*}
\noindent
Since $\boldsymbol{\widetilde{K}}[u^{\varepsilon}]$ is a constant matrix, $\int_{\Omega} \partial_{x_{j}} v \ dx = 0$ for $j \in \{1, 2\}$ and $v \in Y$, it holds:
\begin{align*}
\int_{\Omega} \widetilde{\boldsymbol{K}}[u^{\varepsilon}] \boldsymbol{e}_{1} \cdot \nabla v \ dx =  \sum_{j=1}^{2} (\widetilde{\boldsymbol{K}}[u^{\varepsilon}])_{j1} \int_{\Omega} \partial_{x_{j}} v \ dx = 0.
\end{align*}
To study the fine-to-coarse scale mapping, we need the following related map:
\begin{align}\label{upscaled coefficient mapping}
u \in U \mapsto \boldsymbol{\widetilde{K}}[u] \in \mathbb{R}^{2\times 2}, \quad  u \mapsto   
\begin{pmatrix}
        \boldsymbol{\widetilde{K}}[u]_{11} \ \ \boldsymbol{\widetilde{K}}_{12} \\
        \ \boldsymbol{\widetilde{K}}_{21} \ \ \quad  \boldsymbol{\widetilde{K}}_{22}
\end{pmatrix},
\end{align}
where $\boldsymbol{\widetilde{K}}[u]_{11}= \int_{\Omega} \boldsymbol{K}^{\varepsilon}(\partial_{x_{1}} u + 1 ) \ dx$. Firstly, we prove the following. 
\begin{lemma} \label{lemma, K_tilde} The mapping  $u \in U \mapsto \boldsymbol{\widetilde{K}}[u] \in \mathbb{R}^{2\times 2}$ defined by \eqref{upscaled coefficient mapping} is Lipschitz continuous, i.e., for $u_{1}, u_{2} \in U$ we have
\begin{align} 
\lVert \boldsymbol{\widetilde{K}}[u_{1}] - \boldsymbol{\widetilde{K}}[u_{2}] \rVert & \leq  |\Omega|^{1/2} \beta  \lVert u_{1} - u_{2} \rVert_{U}, \label{Continuity K_tilde} \\
\lVert \boldsymbol{\widetilde{K}}[u_{1}] - \boldsymbol{\widetilde{K}}[u_{2}] \rVert & \leq |\Omega|^{1/2} \varepsilon^{-1} c_{\boldsymbol{K}}  \lVert u_{1} - u_{2} \rVert_{V} \label{Continuity K_tilde in V},
\end{align}
with $\beta>0$ from \eqref{bounds on K}, and the Frobenius norm $\lVert \boldsymbol{\widetilde{K}}[u] \rVert: = \big( \overset{2}{\underset{i,j=1}{\sum}} (\boldsymbol{\widetilde{K}}[u])_{ij}^{2} \big)^{\frac{1}{2}}$ of $\boldsymbol{\widetilde{K}}[u]$.
\end{lemma}
\begin{proof}
The Cauchy--Schwarz inequality and \eqref{bounds on K} yield 
\begin{align}\label{coefficient_cont}
|(\boldsymbol{\widetilde{K}}[u_{1}])_{11}-(\boldsymbol{\widetilde{K}}[u_{2}])_{11}| \leq |\Omega|^{1/2}\beta \lVert \partial_{x_{1}} u_{1} -  \partial_{x_{1}} u_{2}\rVert_{H} \leq |\Omega|^{1/2}\beta \lVert  u_{1} -   u_{2}\rVert_{U} .
\end{align}
In our case, $\boldsymbol{\widetilde{K}}_{12}$, $\boldsymbol{\widetilde{K}}_{21}$, $\boldsymbol{\widetilde{K}}_{22}$ are fixed, hence we get $\lVert \boldsymbol{\widetilde{K}}[u_{1}] - \boldsymbol{\widetilde{K}}[u_{2}] \rVert$ = $| (\boldsymbol{\widetilde{K}}[u_{1}])_{11}  - (\boldsymbol{\widetilde{K}}[u_{2}])_{11}|$ and \eqref{Continuity K_tilde} follows. 

Next, we write $\boldsymbol{\widetilde{K}}[u]_{11}= \int_{\Omega}  \boldsymbol{K}^{\varepsilon}  - (\partial_{x_{1}}\boldsymbol{K}^{\varepsilon})u \ dx$. The Cauchy--Schwarz inequality, $\lVert \nabla \boldsymbol{K}^{\varepsilon} \rVert_{L^{\infty}(\Omega)} \leq \frac{c_{\boldsymbol{K}}}{\varepsilon}$ and the continuous embedding $V \hookrightarrow H$ imply that 
\begin{align*}
\lVert \boldsymbol{\widetilde{K}}[u_{1}] - \boldsymbol{\widetilde{K}}[u_{2}] \rVert  \leq \int_{\Omega} |\partial_{x_{1}}\boldsymbol{K}^{\varepsilon}(u_{2}-u_{1})|dx \leq |\Omega|^{1/2} \varepsilon^{-1} c_{\boldsymbol{K}} \lVert u_{1} - u_{2} \rVert_{V}. 
\end{align*}
\end{proof}
\noindent
We note that \eqref{Continuity K_tilde in V} readily implies \eqref{bi_form_assumption}. Next, we show that Assumption~\ref{Ass:Uniformity} holds, assuming that the influence of $q$ is rather weak. 
\begin{proposition}\label{Upscaled bounds} Suppose that $\lVert q\rVert_{H} \leq \frac{\alpha(\nu + \lVert \nabla u^{\varepsilon} \lVert_{H}^{2})}{ \lVert  u^{\varepsilon}\rVert_{H}}$ for some $0<\nu<1$ and $\alpha \in \mathbb{R}_{+}$ from \eqref{bounds on K}. Then, there exists $B_{\bar{r}}(u^{\varepsilon}) = \{v  \ : \ \lVert u^{\varepsilon} -v \rVert_{U} \leq \bar{r}\} \subset U$ and $C_{b}, C_{c} \in \mathbb{R}_{+}$ such that \eqref{MP: eq 3.1 forms} holds for all $v,w\in Y$ and $u \in B_{\bar{r}}(u^{\varepsilon})$. 
\end{proposition}
\begin{proof} Set $\boldsymbol{v}^{\varepsilon} = (u^{\varepsilon},0)$. Then $\nabla \cdot \boldsymbol{v}^{\varepsilon} = \partial_{x_{1}} u^{\varepsilon}$ and the divergence theorem yield
\begin{align}\label{vanishing_derivative}
\int_{\Omega} \partial_{x_{1}} u^{\varepsilon} \ dx = \int_{\Omega} \nabla \cdot \boldsymbol{v}^{\varepsilon} \  dx = \int_{\partial \Omega} \boldsymbol{v}^{\varepsilon} \cdot \boldsymbol{\eta} \ dx = 0,
\end{align}
since $\restr{u^{\varepsilon}}{\partial \Omega} = 0$. From \eqref{App.A: eq.5}, we deduce
\begin{align}\label{derivation 1}
\boldsymbol{\widetilde{K}}[u^{\varepsilon}]_{11}=\boldsymbol{e}_{1} \cdot \boldsymbol{\widetilde{K}}[u^{\varepsilon}] \boldsymbol{e}_{1}   = \langle \nabla w_{1}^{\varepsilon} \cdot \boldsymbol{K}^{\varepsilon} \nabla w_{1}^{\varepsilon} \rangle_{\Omega} - \langle \nabla u^{\varepsilon} \cdot \boldsymbol{K}^{\varepsilon} \nabla w_{1}^{\varepsilon} \rangle_{\Omega}. 
\end{align}
Integrating the last term in \eqref{derivation 1} by parts gives 
\begin{align}\label{derivation 2}
\int_{\Omega} \nabla u^{\varepsilon} \cdot \boldsymbol{K}^{\varepsilon}\nabla w_{1}^{\varepsilon}  dx  = \int_{\partial \Omega} u^{\varepsilon}  (\boldsymbol{K}^{\varepsilon}\nabla w_{1}^{\varepsilon} \cdot \boldsymbol{\eta})  \ ds - \int_{\Omega} \nabla \cdot (\boldsymbol{K}^{\varepsilon}\nabla w_{1}^{\varepsilon})  u^{\varepsilon} \  dx,
\end{align}
where $\restr{u^{\varepsilon}}{\partial \Omega}=0$ and $-\nabla \cdot (\boldsymbol{K}^{\varepsilon}\nabla w_{1}^{\varepsilon}) = q$ in $\Omega$. Therefore, from \eqref{vanishing_derivative} we obtain 
\begin{align*}
\boldsymbol{\widetilde{K}}[u^{\varepsilon}]_{11} = \langle \nabla w_{1}^{\varepsilon}\cdot \boldsymbol{K}^{\varepsilon}  \nabla w_{1}^{\varepsilon} \rangle_{\Omega}  & - \int_{\Omega}  q u^{\varepsilon}  \ dx   \geq   \alpha \int_{\Omega} |\nabla w_{1}^{\varepsilon}|^{2} dx  - \int_{\Omega}  q u^{\varepsilon}  \ dx  \\ 
 =   \alpha \int_{\Omega} ( |\nabla u^{\varepsilon}|^{2} + 2  \partial_{x_{1}} u^{\varepsilon}  + 1 ) dx  & -  \int_{\Omega}  q u^{\varepsilon} \  dx\geq   \alpha (1 +\lVert \nabla u^{\varepsilon}\rVert_{H}^{2})   - \int_{\Omega}  q u^{\varepsilon} \ dx. 
\end{align*}
The Cauchy--Schwarz inequality and our (amplitude) assumption on $q$ give us
\begin{align*}
\boldsymbol{\widetilde{K}}[u^{\varepsilon}]_{11} \geq   \alpha (1 +\lVert \nabla u^{\varepsilon}\rVert_{H}^{2})  - \lVert q \rVert_{H} \lVert u^{\varepsilon} \rVert_{H}  \geq (1-\nu)\alpha=:C_{\alpha}.
\end{align*}
Integration by parts, using the equation satisfied by $u^{\varepsilon}$, the Cauchy--Schwarz inequality, and \eqref{grad_bound} result in the estimate
\begin{align*}
|\boldsymbol{\widetilde{K}}[u^{\varepsilon}]_{11} \big| \leq  \beta  +   \beta \lVert  \nabla u^{\varepsilon} \rVert_{H}  \leq \beta +   \beta (\frac{\beta + c_{p} \lVert q \rVert_{H}}{\alpha})=: C_{\beta}.
\end{align*}
Let $0 \leq \bar{s} \leq \nu(1-\nu)\alpha $ and $\bar{r}=\frac{\bar{s}}{|\Omega|^{1/2} \beta}$. Then \eqref{Continuity K_tilde} implies that for $u \in B_{\bar{r}}(u^{\varepsilon})$ it holds:
\begin{align}\label{sequence_bounds}
0 < (1-\nu) C_{\alpha} \leq \boldsymbol{\widetilde{K}}[u^{\varepsilon}]_{ii} - \bar{r} \leq \boldsymbol{\widetilde{K}}[u]_{11} \leq \boldsymbol{\widetilde{K}}[u^{\varepsilon}]_{11} +  \bar{r} \leq C_{\beta} + \nu(1-\nu)\alpha. 
\end{align}
Then, we set $C_{c}:=(1-\nu)C_{\alpha}$ and $C_{b}:=C_{\beta} + \nu(1-\nu)\alpha$ and apply the standard coercivity and continuity estimates to complete the proof.
\end{proof}
\noindent
The Lax--Milgram lemma and Proposition~\ref{Upscaled bounds} imply that for $u \in B_{\bar{r}}(u^{\varepsilon})$ there exist a unique solution $y(u) \in Y$ of \eqref{weak formulation} with $\lVert y(u) \rVert_{Y} \leq C\lVert q \Vert_{H}$, where $C \in \mathbb{R}_{+}$ is independent of $u$, but also of $\varepsilon$, since $\bar{r}$, $C_{c}$ and $C_{b}$ are independent of $\varepsilon$ in Proposition~\ref{Upscaled bounds}. The constant coefficient in \eqref{weak formulation} implies that $y(u) \in U$ for all $u \in B_{\bar{r}}(u^{\varepsilon})$, hence we use \eqref{MP: eq 8}. Proposition~\ref{Upscaled bounds} motivates the following rather general assumption.
\begin{Assumption}
\label{Ass: Data for uniformity}
For $q \in H$, $\boldsymbol{K}^{\varepsilon}\in C^{0, 1}(\bar{\Omega})$, one can find $\bar{r}^{\varepsilon} \in \mathbb{R}_{+}$ and respective $B_{\bar{r}^{\varepsilon}}(u^{\varepsilon}) \subset U$ such that $y(u) \in Y$ exists and $\lVert y(u) \rVert_{Y} \leq C \lVert q \rVert_{H}$ for all $u \in B_{\bar{r}^{\varepsilon}}(u^{\varepsilon})$, where $C \in \mathbb{R}_{+}$ is independent of $u$. 
\end{Assumption}
Next, we verify Assumption~\ref{Ass: Continuity} about the continuity of our fine-to-coarse scale map.  
\begin{proposition}
\label{Th: continuity sol map}
Suppose that Assumption~\ref{Ass: Data for uniformity} holds and $\{u_{k}^{\varepsilon}\} \subset B_{\bar{r}^{\varepsilon}}(u^{\varepsilon}) \cap X$ is an approximating sequence of $u_{\varepsilon}$. Then $S(u_{k}^{\varepsilon}) \rightharpoonup S(u^{\varepsilon})$ in $Y$ as $k \rightarrow \infty$, where $S : U \rightarrow Y$  is the fine-to-coarse scale map.
\end{proposition}
\begin{proof}
Note that $\lVert y(u_{k}^{\varepsilon}) \rVert_{Y} \leq C \lVert q \rVert_{H}$ for $\{u_{k}^{\varepsilon}\}$, $C \in  \mathbb{R}_{+}$ and $k \in \mathbb{N}$ as stated in Assumption~\ref{Ass: Data for uniformity}. The reflexivity of $Y$ and the Banach–-Alaoglu theorem imply that $\{y(u_{k}^{\varepsilon}) \}$ admits a weakly convergent subsequence with its elements still denoted by $y(u_{k}^{\varepsilon})$. Let $\hat{y} \in Y$ denote that weak limit. Rearranging terms in \eqref{weak formulation}, we obtain
\begin{align*}
b_{\mathcal{L}}[u^{\varepsilon}](y(u^{\varepsilon}),v) + \int_{\Omega} \big( \boldsymbol{\widetilde{K}}[u_{k}^{\varepsilon}] - \boldsymbol{\widetilde{K}}[u^{\varepsilon}] \big)\nabla y(u^{\varepsilon}) \cdot \nabla v \ dx =  \int_{\Omega}qv \ dx, \quad \forall{v} \in Y.
\end{align*}
Note that $b_{\mathcal{L}}[u^{\varepsilon}](y(u_{k}^{\varepsilon}),v) \rightarrow b_{\mathcal{L}}[u^{\varepsilon}](\hat{y},v)$ as $k \rightarrow \infty$ for all $v \in Y$, since $\mathcal{L}[u^{\varepsilon}]\in L(Y,Y^{*})$ and hence it is weakly continuous. The continuity \eqref{Continuity K_tilde} and the Cauchy--Schwarz inequality yield
\begin{align*}
\int_{\Omega} \big( \boldsymbol{\widetilde{K}}[u_{k}^{\varepsilon}] - \boldsymbol{\widetilde{K}}[u^{\varepsilon}] \big)\nabla y(u_{k}^{\varepsilon}) \cdot \nabla v \ dx \leq C\lVert u_{k}^{\varepsilon} - u^{\varepsilon} \rVert_{U} \lVert q \rVert_{H} \lVert v \rVert_{Y} \rightarrow 0 \ \text{as} \  k\rightarrow \infty,
\end{align*}
since $C\in \mathbb{R}_{+}$ does not depend on $u_{k}^{\varepsilon}$ and $u^{\varepsilon}$. This shows that $\hat{y}=y(u^{\varepsilon})$. 
\end{proof}

The continuity of $S^{\prime} : B_{\bar{r}^{\varepsilon}}(u^{\varepsilon}) \subset U \rightarrow L(U, Y)$ can be established through standard techniques, which we briefly outline here by examining the sensitivity equation. The sensitivity $z:=S^{\prime}(u)h \in Y$ of $S(\cdot)$ at $u\in B_{\bar{r}^{\varepsilon}}(u^{\varepsilon})$ in the direction $h \in U$ is given as the solution of the linearized state equation 
\begin{align}\label{Sens WF}
\langle e_{y}(y(u), u) z, v \rangle_{Y^{\ast},Y} = - \langle e_{u}(y(u), u)h, v \rangle_{Y^{*}, Y}  \quad \forall{v} \in Y,
\end{align}
where the partial derivatives $e_{y}(y, u): Y \rightarrow Y^{\ast}$ and $e_{u}(y, u): U \rightarrow Y^{\ast}$ read  
\begin{align*}
\langle e_{y}(y, u)w, v \rangle_{Y^{\ast},Y}  = \int_{\Omega} \widetilde{\boldsymbol{K}}[u] \nabla w \cdot \nabla v \ dx, \ \langle e_{u}(y, u)h, v \rangle_{Y^{\ast}, Y}  = \int_{\Omega} \widetilde{\boldsymbol{K}}_{u}[h] \nabla y \cdot \nabla v \ dx,
\end{align*}
respectively, and $\widetilde{\boldsymbol{K}}_{u}[h] \in \mathbb{R}^{2\times 2}$ is given by $(\widetilde{\boldsymbol{K}}_{u}[h])_{11} = \int_{\Omega} \boldsymbol{K}^{\varepsilon}\partial_{x_{1}} h \ dx$ and $(\widetilde{\boldsymbol{K}}_{u}[h])_{21}=(\widetilde{\boldsymbol{K}}_{u}[h])_{12}=(\widetilde{\boldsymbol{K}}_{u}[h])_{22}=0$ with $\lVert \widetilde{\boldsymbol{K}}_{u}[h] \rVert  \leq  |\Omega|^{1/2} \beta  \lVert h \rVert_{U}$ for $h \in U$. Here, Proposition~\ref{Upscaled bounds} or Assumption~\ref{Ass: Data for uniformity} implies that $e_{y}(y, u)$ is boundedly invertible for all $u\in B_{\bar{r}^{\varepsilon}}(u^{\varepsilon})$, verifying Assumption~\ref{Ass: System invertibility} for the self-adjoint operator $\mathcal{L}[u]$. For $ u_{1}, u_{2} \in B_{\bar{r}^{\varepsilon}}(u^{\varepsilon})$, let $z_{1}= S'(u_{1})h$ and $z_{2}= S'(u_{2})h$. Invoking the well-posedness of \eqref{Sens WF} and \eqref{Continuity K_tilde} while estimating the right-hand side of \eqref{Sens WF} provide us with the estimate 
\begin{align*}
\lVert \big( S'(u_{1}) - S'(u_{2})\big) h \rVert_{Y} \leq  C_{S'} \lVert u_{1} - u_{2} \lVert_{U} \lVert h \rVert_{U},
\end{align*}
where $C_{S'} \in \mathbb{R}_{+}$ generally depends on $u$. Therefore, $S(u)$ is continuously Fr\'echet differentiable and Assumption~\ref{Ass: differentiability} holds true.

\subsection{Implementation issues}
We provide the implementation details of Algorithm~\ref{alg:hybrid learning} for problem \eqref{HC eq 1} with a periodic coefficient $\boldsymbol{K}^{\varepsilon}(x)$. Firstly, the derivative $e_{\boldsymbol{\theta}}(y, \boldsymbol{\theta}) \in L(\mathbb{R}^{N_{n}},Y^{*})$ is given by
\begin{align*}
\langle e_{\boldsymbol{\theta}}(y, \boldsymbol{\theta})\boldsymbol{s}, v \rangle_{Y^{*}, Y}  & = \int_{\Omega} \widetilde{\boldsymbol{K}}_{\boldsymbol{\theta}}[\boldsymbol{s}] \nabla y \cdot \nabla v \ dx, \quad  \widetilde{\boldsymbol{K}}_{\boldsymbol{\theta}}[\boldsymbol{s}] =  \int_{\Omega} \boldsymbol{K}^{\varepsilon} \langle \nabla_{\boldsymbol{\theta}}   (\partial_{x_{1}} v_{\boldsymbol{\theta}, n}), \boldsymbol{s} \rangle_{\mathbb{R}^{N_{n}}} \ dx
\end{align*}
for $\boldsymbol{s}\in \mathbb{R}^{N_{n}}$.
Given the $k$-th unit vector $\boldsymbol{e}_{k} \in \mathbb{R}^{N_{n}}$, one obtains
\begin{align} \label{adj_contr}
\langle y'(\boldsymbol{\theta})^{\ast} \partial_{y} J(y(\boldsymbol{\theta}),\boldsymbol{\theta}), \boldsymbol{e}_{k}  \rangle_{\mathbb{R}^{N_{n}}} =  \langle e_{\boldsymbol{\theta}}(y(\boldsymbol{\theta}), \boldsymbol{\theta}) \boldsymbol{e}_{k}, p\rangle_{Y^{\ast}, Y}.
\end{align}
The formula \eqref{adj_contr}  is useful for the assembly of the first summand in \eqref{DOCP: eq 5}.  The algebraic FE systems are described in Section~\ref{Section 2}, but we note that $\mathbb{B}_{h}[\boldsymbol{\theta}] = \widetilde{\boldsymbol{K}}[v_{\boldsymbol{\theta}, n}] \mathbb{A}_{h}$ with $(\mathbb{A}_{h})_{ij}: =  \langle  \nabla \phi_{i}, \nabla \phi_{j} \rangle_{H}$. The discrete counterpart of \eqref{adj_contr} is given by
\begin{align*}
\langle e_{\boldsymbol{\theta}}(y_{h}(\boldsymbol{\theta}), \boldsymbol{\theta})\boldsymbol{e}_{k}, p_{h}\rangle_{Y^{\ast}, Y} = \boldsymbol{y}_{h}^{T} \mathbb{E}_{h}[\boldsymbol{\theta}_{k}] \ \boldsymbol{p}_{h}, \quad 1 \leq k \leq N_{n}, 
\end{align*}
where $\mathbb{E}_{h}[\boldsymbol{\theta}_{k}] \in \mathbb{R}^{N_{h}\times N_{h}}$, $(\mathbb{E}_{h}[\boldsymbol{\theta}_{k}])_{ij}:=\langle e_{\boldsymbol{\theta}}(\phi_{i}, \boldsymbol{\theta})\boldsymbol{e}_{k}, \phi_{j} \rangle_{Y^{\ast}, Y}$. We get $(\mathbb{E}_{h}[\boldsymbol{\theta}])_{ij}   = \widetilde{\boldsymbol{k}}_{M}[\boldsymbol{\theta}] (\mathbb{A}_{h})_{ij}$, where $\widetilde{\boldsymbol{k}}_{M}[\boldsymbol{\theta}] \in \mathbb{R}^{N_{n}}$ represents approximations of $\widetilde{\boldsymbol{K}}_{\boldsymbol{\theta}}[\boldsymbol{e}_{k}]$ by the Monte-Carlo approach and using the residual collocation points 
\begin{align*}
(\widetilde{\boldsymbol{k}}_{M}[\boldsymbol{\theta}])_{j}= \frac{1}{M_{\Omega}}\sum_{i=1}^{M_{\Omega}}\boldsymbol{K}^{\varepsilon}(x_{1, i}^{r}, x_{2, i}^{r}) \frac{\partial^{2}v_{\boldsymbol{\theta}, n}^{M}(x_{1, i}^{r}, x_{2, i}^{r})}{\partial \boldsymbol{\theta}_{j} \ \partial x_{1}}, \quad 1\leq j \leq n. 
\end{align*}

We also use a multiscale Fourier feature network \cite{tancik2020fourier, wang2021eigenvector} to reduce the effect of spectral bias. This architecture includes Fourier feature mappings $\mathcal{F}^{(k)}: \mathbb{R} \rightarrow \mathbb{R}^{2m}$:
\begin{align*}
\mathcal{F}^{(k)}(x) = \begin{pmatrix}
\cos (2 \pi \boldsymbol{B}^{(k)}x), \
\sin (2 \pi \boldsymbol{B}^{(k)}x)
\end{pmatrix}, \quad 1 \leq k \leq K,
\end{align*}
where each entry of $\boldsymbol{B}^{(k)} \in \mathbb{R}^{m \times d}$ is sampled from a Gaussian distribution $\mathcal{N}(0, \varrho_{k}^{2})$ with $\varrho_{k} > 0$ a specified hyperparameter. These features are used as inputs for the hidden layers, which are defined for $ 1\leq k \leq K$ and $ 2 \leq l \leq L-1$ as follows:
\begin{align*}
z_{1}^{(k)}  = \sigma( W_{1} \mathcal{F}^{(k)}(x) + b_{1}), \
z_{l}^{(k)}  = \sigma( W_{l} z_{l-1}^{(k)} + b_{l}). 
\end{align*}
Next, the above outputs are concatenated within the linear layer as follows:
\begin{align*}
v_{\boldsymbol{\theta}, n}^{M} = W_{L} \big[z_{L}^{(1)},...,z_{L}^{(K)}\big] + b_{L},
\end{align*}
where $W_{L}$ and $b_{L}$ are the weights and biases of the output layer. To enforce the boundary conditions exactly \cite{lagari2020systematic, lu2021physics, sukumar2022exact}, one may set $v_{\boldsymbol{\theta}, n}^{l,M} := l(x) v_{\boldsymbol{\theta}, n}^{M}$, where $l$ is the corresponding signed distance function for $\partial \Omega$.

\subsection{Numerical results} For our numerical experiments, we use the following setup unless otherwise specified. The multiscale Fourier feature network (Ms-PINN) is used as the main architecture of choice, with the two Fourier features initialized by $\varrho_{1} = 1$ and $\varrho_{2} = 1/\varepsilon$, two hidden layers with 128 neurons each, and $\tanh{(x)}$ as the activation functions. The full batch is used for training with the Adam algorithm \cite{kingma2014adam}, and  $\beta_{1}^{Ad}=0.9$ and $\beta_{2}^{Ad}=0.999$ are chosen as the hyperparameters. The learning rate for the exponential learning rate schedule is initialized as $lr(0)=5e-4$, with a decay-rate of 0.75 every 1000 training iterations. For PINNs and our hybrid approach, we use automatic loss balancing \cite[Algorithm 1(c)]{wang2023expert} at every 100 steps. We choose \eqref{MP: eq 8} for coupling, the uniform quadrature rule for \eqref{full_discrete_couple_term}, and $\delta = \varepsilon$ for \eqref{full_discrete_comp}. To reduce computational costs, we update the coarse-scale information only if $|\widetilde{\boldsymbol{K}}[v_{\boldsymbol{\theta}^{(it+T_{K})},n}^{M}] - \widetilde{\boldsymbol{K}}[v_{\boldsymbol{\theta}^{(it)},n}^{M}]| > 10^{-2}$, where $T_{K}=100$. Our implementation is based on JAX \cite{bradbury2018jax} and FEniCS Dolphin \cite{logg2010dolfin}. 

The following example is used to illustrate the concept of upscaling consistency: 
\begin{align}\label{1D problem numerics}
-&\partial_{x}(\boldsymbol{K}^{\varepsilon}  \partial_{x}u^{\varepsilon}) = f^{\varepsilon},  \quad  \text{in} \  \Omega:=(0,1),  \quad \text{and} \quad   u^{\varepsilon}(0)  = 0  ,   \   \ u^{\varepsilon}(1)  = 0, 
\end{align}
where $\boldsymbol{K}^{\varepsilon}(x)=2 + \sin(2\pi \varepsilon^{-1} x)$, $f^{\varepsilon}: = q + \partial_{x} \boldsymbol{K}^{\varepsilon}$ and $q = -3(2x - 1)$. Then, $\widetilde{\boldsymbol{K}}[u] = \int_{0}^{1} \boldsymbol{K}^{\varepsilon} (\partial_{x} u +1 ) \ dx$ and $\boldsymbol{K}^{\ast} = (\int_{0}^{1} \frac{dx}{\boldsymbol{K}^{\varepsilon}(x)})^{-1}=\sqrt{3}\approx 1.732$, see \cite[Section 1.2]{jikov2012homogenization} for the latter. For $\varepsilon=1/15$ and $\varepsilon=0.01$, we get $\widetilde{\boldsymbol{K}}[u^{\varepsilon}_{h}]\approx 1.733$ and $\widetilde{\boldsymbol{K}}[u^{\varepsilon}_{h}] \approx 1.732$, where $u^{\varepsilon}_{h}$ are computed using finite elements on equidistant meshes with $h = 5 \cdot 10^{-4}$ and $h = 2 \cdot 10^{-4}$, respectively (these meshes are further used for validation). Then, the influence caused by an upscaling method in Assumption~\ref{Ass: Consistent parametrization} is negligible in the context of Theorem~\ref{Th: Upscaling consistency}. Setting $it_{\text{max}} = 2 \cdot 10^{5}$, we use $M_{\Omega}=2000$ and $M_{\Omega}=5000$ uniformly random collocation points for $\varepsilon=1/15$ and $\varepsilon=0.01$ (both with hard boundary constraints), respectively, and the mesh size $h=0.004$ for the coarse-scale problem.  We use $l(x)=4x(1-x)$ for hard boundary constraints imposition. 

\begin{figure}[H]
    \centering
    \includegraphics[width=1\textwidth]{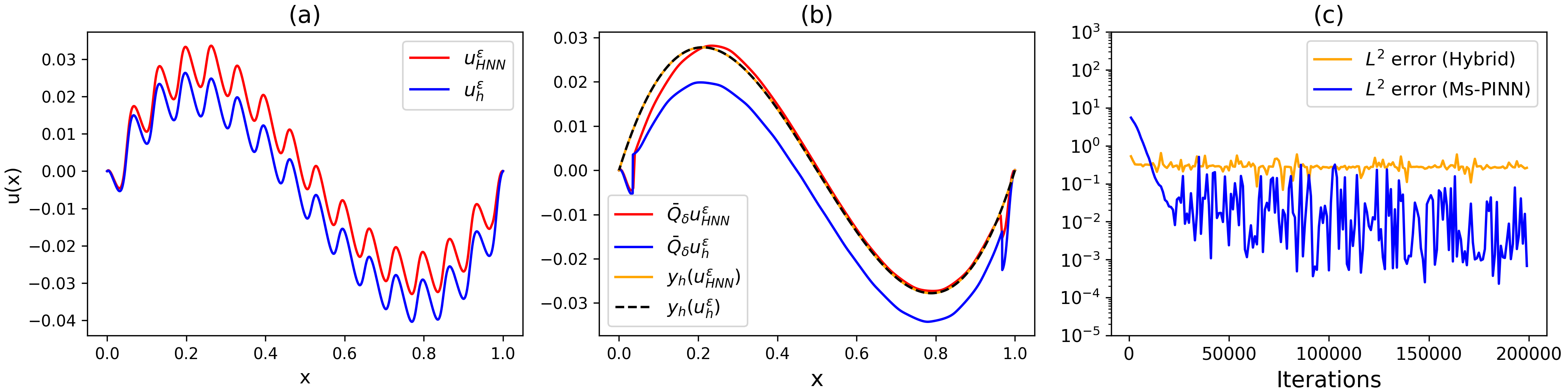}
    \vspace*{-6.7mm}
    \caption{\small{Results for $\varepsilon = 1/15$, 1D problem: \textbf{(a)} Hybrid fine-scale solution $u_{HNN}^{\varepsilon}$ and FEM solution $u_{h}^{\varepsilon}$, \textbf{(b)} Predicted state $y_{h}(u_{HNN}^{\varepsilon})$, true discrete state $y_{h}(u_{h}^{\varepsilon})$, compression $\bar{Q}_{\delta}u_{HNN}^{\varepsilon}$ and $\bar{Q}_{\delta}u_{h}^{\varepsilon}$, \textbf{\textbf{c}} Relative $L^{2}$ error of the hybrid solver and Ms-PINN vs iterations.}}
    \label{fig: sim 1}
\end{figure}

If $\varepsilon$ is not small enough, one expects an upscaling consistency gap between $y_{h}(u)$ and $\bar{Q}_{\delta}u$ for $u=u_{h}^{\varepsilon}$, and this gap is typically overfitted by a neural network $u=v^{M}_{\boldsymbol{\theta},n}$. This overfitting can introduce discrepancies in the PDE residual of \eqref{COCP: eq 1}. The latter can impede the route to convergence on a fine scale, rendering it generally not possible. Indeed, for $\varepsilon = 1/15$ in \eqref{1D problem numerics}, we observe in Fig~\ref{fig: sim 1}(c) that the relative $L^{2}$ error of the hybrid solver oscillates around $0.3$ on our validation set, and does not drop significantly below that value, compared to the standard Ms-PINN approach without coarse-scale constraints, which shows $\sim 10^{-2}$ accuracy. We can see in Fig~\ref{fig: sim 1}(b) that $y_{h}(v^{M}_{\boldsymbol{\theta},n})$ matches $y_{h}(u_{h}^{\varepsilon})$ well, but $\bar{Q}_{\delta}v^{M}_{\boldsymbol{\theta},n}$ is much closer to $y_{h}(v^{M}_{\boldsymbol{\theta},n})$ compared to $\bar{Q}_{\delta}u_{h}^{\varepsilon}$, thus violating the upscaling consistency and causing a significant approximation error in Fig~\ref{fig: sim 1}(a).     
\begin{figure}[H]
    \centering
    \includegraphics[width=1\textwidth]{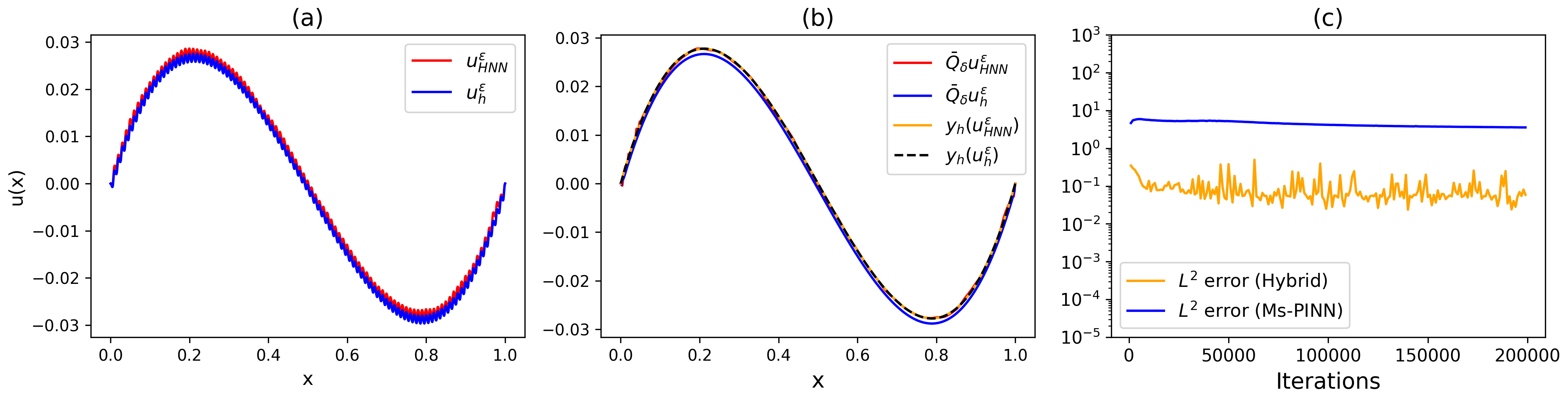}
    \vspace*{-6.7mm}
    \caption{\small{Results for $\varepsilon = 0.01$, 1D problem: same description as for Fig~\ref{fig: sim 1}}.}
    \label{fig: sim 2}
    \vspace*{-3mm}
\end{figure} 

Theorem~\ref{Th: Upscaling consistency} shows that the upscaling consistency gap decreases with decreasing $\varepsilon$, thus softening the conflict between the PDE residual term and the coupling term in the optimization and reducing the effect of coarse-scale overfitting. Indeed, we observe a smaller gap between $y_{h}(u^{\varepsilon}_{h})$ and $\bar{Q}_{\delta}u^{\varepsilon}_{h}$ for $\varepsilon=0.01$ in Fig~\ref{fig: sim 2}(b), and a lower resulting $L^{2}$ error for the hybrid method in Fig~\ref{fig: sim 2}(c), compared to the previous example. In addition, the hybrid solver outperforms the Ms-PINN approach by almost two orders of magnitude in terms of relative error. However, the hybrid solver error still oscillates around $0.1$: there is a visible overfitting of the coarse-scale component at $x\approx 0.2$ and $x\approx 0.8$ in Fig~\ref{fig: sim 2}(b), which causes visible approximation errors in Fig~\ref{fig: sim 2}(a) at the same spatial points for the fine-scale hybrid approximation. 
\begin{figure}[H]
    \centering
    \includegraphics[width=1\textwidth]{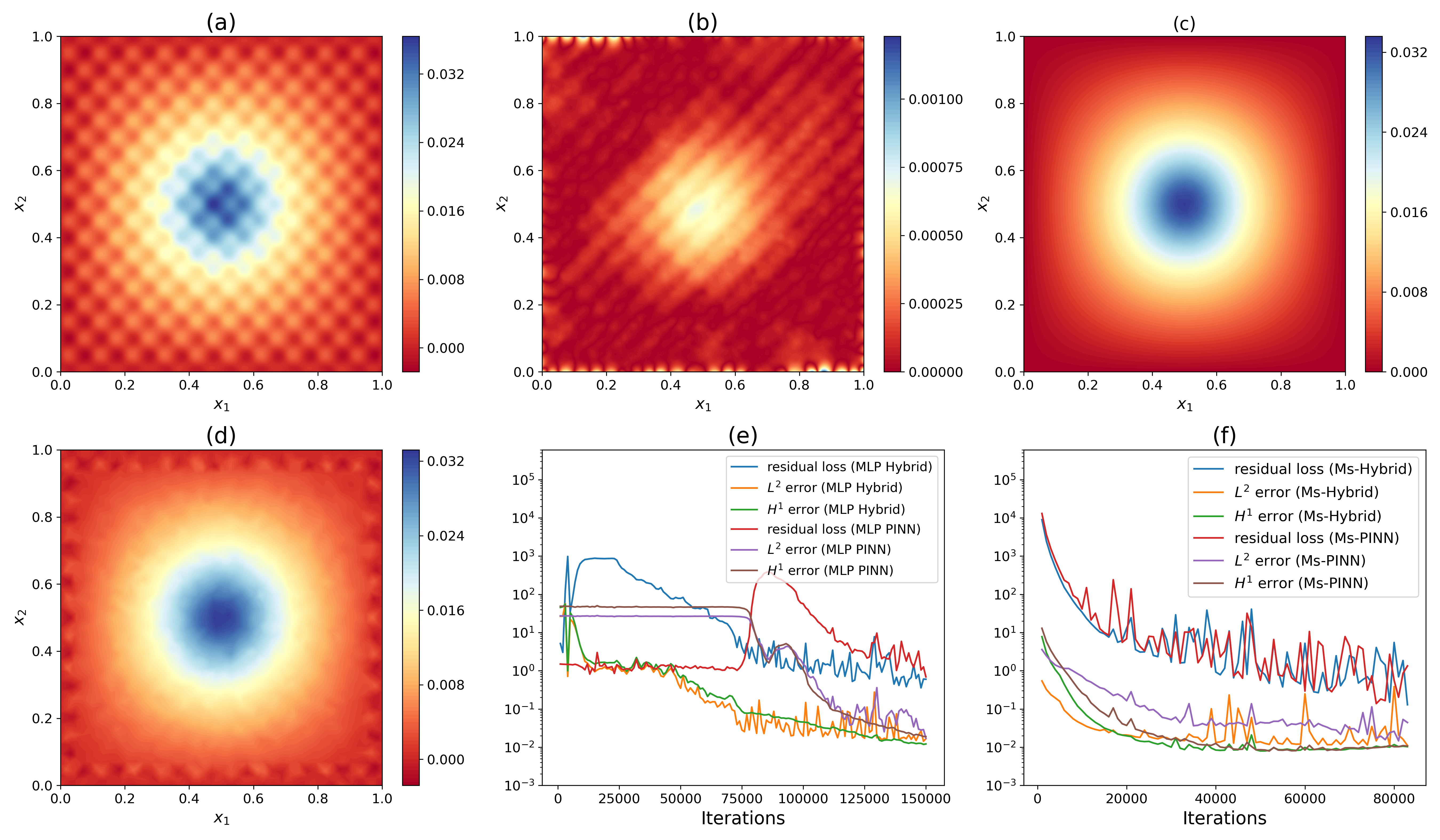}
    \vspace*{-6.5mm}
    \caption{\small{Results for $\varepsilon = 0.1$, 2D problem: \textbf{(a)} Hybrid MLP-based fine-scale solution $u_{HNN}^{\varepsilon}$, \textbf{(b)} Pointwise error $|u_{HNN}^{\varepsilon}(x) - u_{h}^{\varepsilon}(x)|$,  \textbf{(c)} Predicted state $y_{h}(u_{HNN}^{\varepsilon})$, \textbf{(d)} Compression $\bar{Q}_{\delta}u_{HNN}^{\varepsilon}$, \textbf{(e)} The PDE residual losses and relative $L^{2}$ and $H^{1}$ errors vs iterations for the MLP-based hybrid solver and the MLP PINN, \textbf{(f)} The PDE residual losses and relative $L^{2}$ and $H^{1}$ errors vs iterations for the Ms-PINN-based hybrid solver and the Ms-PINN.}}
    \label{fig: sim 3}
\end{figure}

Next, we study the 2D problem \eqref{HC eq 1} with the isotropic coefficient $\boldsymbol{K}^{\varepsilon}(x) = 3 + \sin(\pi \varepsilon^{-1} x_{1}) + \sin (2 \pi \varepsilon^{-1} x_{2})$ with $\varepsilon = 0.1$ and  $q(x)=10\exp{(-100\lVert x - 0.5\rVert_{2}^{2})}$. We note that $\boldsymbol{K}^{\varepsilon}$ is periodic, hence we need to determine $\boldsymbol{\widetilde{K}}_{11}$ only. In addition to the Ms-PINN, we employ the multilayer perceptron (MLP) PINN with 8 hidden layers with 64 neurons each for both the standard and hybrid approaches while setting $lr(0)=1e-3$ for the initial learning rate. We use $M_{\Omega}=25000$ and $M_{\partial \Omega}$=2500 uniformly random collocation points and the mesh size $h\approx 0.04$ for the coarse-scale discretization. Figure~\ref{fig: sim 3} shows exemplary realizations of the hybrid solver. The MLP network, prone to low-frequency bias, benefits from embedding coarse-scale constraints, which accelerates the approximation of the coarse-scale component and mitigates the bias, as seen in Fig~\ref{fig: sim 3}(e). The MLP approach struggles with being trapped at poor local minima, failing to provide a reasonable approximation for around 75k iterations. In contrast, the hybrid approach exhibits steadily decreasing residuals and errors. For the Ms-PINNs, the hybrid method also shows faster decay of the $L^{2}$ and $H^{1}$ relative errors\footnote{The $H^{1}$ norm is induced by the $L^{2}$ norm of the gradient.}, though the differences are less pronounced compared to the MLP approach, as the Ms-PINN approach is less prone to the spectral bias.
\begin{table}
\caption{Comparison of the Ms-PINNs and the MLP-based PINNs with the hybrid method  for the 2-D problem. The presented quantities are averages over 5 random initializations of NN parameters and collocation points. For the MLP PINN only one  realization is presented, as all others failed to deliver reasonable results. All runs are executed on A100 GPUs with 40 GB RAM.}\label{tab:2} 
\vspace*{-1mm}
\begin{center}[H]
\begin{tabular}{|  l | l | l | l | l |}
\hline
Method &  Rel. $L^{2}$ error  & Rel. $H^{1}$ error  & iterations
& time (sec) \\ \hline
Hybrid Ms-PINN  &  1.28E-02 $\pm$ 0.003 &  2.23E-02 $\pm$  0.008 & 77478 & 554 \\
Hybrid MLP  &  1.42E-02 $\pm$ 0.001
 & 1.44E-02
 $\pm$ 0.008 & 139400 & 1725
 \\ 
Ms-PINN  &  1.97E-02 $\pm$ 0.005 &  2.20E-02 $\pm$ 0.005 & 70009 & 354 \\ 
MLP PINN  & 1.71E-02 $\pm$ 1.000
 &  1.89E-02 $\pm$ 1.000
 &  150000
 &  1221 \\ 
\hline
\end{tabular}
\end{center} 
\vspace*{-4.5mm}
\end{table}

In Table 1, we report the $L^{2}$ and $H^{1}$ relative errors corresponding to the smallest residual norm value within the learning history, the number of iterations required to achieve it, and the compute time, all averaged over 5 random initializations. According to our selection criteria, the hybrid Ms-PINN requires slightly more iterations than the standard Ms-PINN, but in Fig~\ref{fig: sim 3}(f) the convergence to the reported values in Table 1 is factually noticeably faster for the hybrid approach, suggesting more robustness at the cost of compute time. For the standard MLP PINN, one successful realization was obtained for Fig~\ref{fig: sim 3}(e), while many others failed to approximate the solution. The hybrid MLP network consistently yields robust results across all realizations.

\section{Conclusion}

This paper focuses on the structural properties of a learning-informed PDE-constrained optimization problem with a PINN-based objective, resulting in a hybrid two-scale solver. Our approach integrates conventional multiscale methodologies and deep learning techniques, and the proposed PDE-constrained optimization setting seems particularly well-suited for this purpose. However, selecting a suitable neural network approximation scheme and developing an efficient optimization algorithm, aimed at enhancing accuracy while minimizing computational time, are both essential for further algorithmic developments. This task requires considering recent advancements in PINNs and the field of PDE-constrained optimization. For example, one may consider using non-smooth activation functions, introducing non-smoothness into the related PDE-constrained optimization. Such a setting challenges both the derivation of suitable stationarity conditions for the PDE-constrained multiscale approach and its numerical solution. A comprehensive investigation of such a setting, as well as the development of efficient optimization algorithms to address other problems, remain part of our future work in this area.

\bibliographystyle{siamplain}
\bibliography{references}
\end{document}